\documentclass[11pt]{article}

\usepackage{amsmath,amssymb,amsthm}
\usepackage{mathrsfs}
\usepackage{bbm}
\usepackage{amsfonts}
\usepackage{xcolor}

\usepackage{mathtools}
\usepackage{xr-hyper}
\newcommand{\R}{ \mathbb{R}}          % set of real numbers
        
\newcommand{\C}{ \mathbb{C}}          % set of complex numbers
\newcommand{\N}{ \mathbb{N}}          % set of natural numbers
\newcommand{\Z}{ \mathbb{Z}}          % set of integers
\def\V#1{{\boldsymbol{#1}}}           % vector-valued variables
\def\M#1{{\bf{#1}}}                   % matrices or fixed vectors
\newcommand{\Tr}{\mathsf{T}}          % transpose
\newcommand{\Fourier}{\Spc F}         % Fourier transform

\def\Spc#1{{\mathcal{#1}}}            % function spaces
\def\Op#1{{\mathrm{#1}}}              % operator
\def\ee{\mathrm{e}}                   % exponential
\def\jj{\mathrm{i}}                   % square root of -1
\newcommand{\dint}{{\rm d}}           % integration element

\newcommand{\embedD}{\overset{d}{\hookrightarrow}} % continuous & dense embedding
\newcommand{\embedC}{\overset{}{\hookrightarrow}} % continuous embedding
\newcommand{\eqdef}{:=}
\newcommand{\Lop}{\Op L}
\newcommand{\Dop}{\Op D}
\newcommand{\Identity}{\mathrm{Id}}
\newcommand{\bw}{\V \omega}

\providecommand{\Prob}{}\renewcommand{\Prob}{\mathscr{P}}
\def\ProbF{\widehat{\Prob}} 
\def\Exp{\mathop{\mathbb{E}}\nolimits}% expectation
\def\wRand{\xi}
\def\bwRand{{\boldsymbol \xi}}
\def\Rand#1{{\uppercase{#1}}}         % random variables
\def\VRand#1{{\uppercase{\boldsymbol #1}}}   % random vectors
\def\gRand#1{{\uppercase{ #1}}}       % generalized random processes
\theoremstyle{plain}
\newtheorem{theorem}{Theorem}[section]

\newtheorem{corollary}[theorem]{Corollary}
\newtheorem{proposition}[theorem]{Proposition}
\theoremstyle{definition}
\newtheorem{definition}[theorem]{Definition}

\theoremstyle{remark}
\newtheorem{remark}[theorem]{Remark}

\usepackage[numbers,sort&compress]{natbib}
\usepackage[hidelinks]{hyperref}
\usepackage{clchanges}   
\begin{document}

\title{Variational Splines and Generalized Gaussian Processes}

%\title{Splines and Operators:\\ Part III---
%Generalized Gaussian Processes\\ and Bayesian Estimation}
\author{Michael Unser\thanks{Biomedical Imaging Group, \'Ecole polytechnique
f\'ed\'erale de Lausanne (EPFL), CH-1015 Lausanne, Switzerland
(michael.unser@epfl.ch).}}
\date{}
\maketitle

\begin{abstract}
Tikhonov regularization, kriging, and Gaussian-process regression share the same algorithm
because, for a linear inverse problem with a Gaussian prior and Gaussian noise, the
minimum-mean-square-error (MMSE) estimator is the solution of a regularized least-squares
problem. We show that this equivalence survives the passage from vectors to generalized
functions: the regularized least-squares estimators become generalized splines, and the Gaussian
random vectors become generalized Gaussian processes (GGPs) on the dual of a nuclear space. A
whitening operator $\Lop: \Spc S \to L_2$ specifies both the native Hilbert space
$\Spc H_\Lop \subset \Spc S'$ of the splines and a family of GGPs, which solve a stochastic
differential equation subject to boundary conditions. We prove that the Gaussian spaces of these
processes are congruent to the core $\Lop^\ast(L_2)$ of the native space, and that their MMSE
estimators from finitely many linear measurements are generalized splines.
Our framework extends the classical correspondences of Kimeldorf--Wahba, Parzen and Kailath in
two respects: our estimators form a broader class that includes generalized functions, and we
give a proper Bayesian treatment of the case where the null space of $\Lop$ is non-trivial.
The framework gives a precise meaning to statements such as ``fractional splines are the
optimal estimators of fractional Brownian motion,'' and it provides a rigorous Bayesian
formulation of infinite-dimensional linear inverse problems.

%% Abstract for ArXiv
%Tikhonov regularization, kriging, and Gaussian-process regression share the same algorithm
%because, for a linear inverse problem with a Gaussian prior and Gaussian noise, the
%minimum-mean-square-error (MMSE) estimator is the solution of a regularized least-squares
%problem. We show that this equivalence survives the passage from vectors to generalized
%functions: the regularized least-squares estimators become generalized splines, and the Gaussian
%random vectors become generalized Gaussian processes (GGPs) on the dual of a nuclear space. A
%whitening operator $\mathrm{L}: \mathcal{S} \to L_2$ specifies both the native Hilbert space
%$\mathcal{H}_\mathrm{L} \subset \mathcal{S}'$ of the splines and a family of GGPs, which solve a stochastic
%differential equation subject to boundary conditions. We prove that the Gaussian spaces of these
%processes are congruent to the core $\mathrm{L}^\ast(L_2)$ of the native space, and that their MMSE
%estimators from finitely many linear measurements are generalized splines.
%Our framework extends the classical correspondences of Kimeldorf-Wahba, Parzen and Kailath in
%two respects: our estimators form a broader class that includes generalized functions, and we
%give a proper Bayesian treatment of the case where the null space of $\mathrm{L}$ is non-trivial.
%The framework gives a precise meaning to statements such as ``fractional splines are the
%optimal estimators of fractional Brownian motion,'' and it provides a rigorous Bayesian
%formulation of infinite-dimensional linear inverse problems.

\end{abstract}
\pagebreak
\tableofcontents
\bigskip

\section{Introduction}
\label{Sec:Introduction}
The equivalence between optimal estimation under the Gaussian hypothesis and regularized least
squares is one of the organizing principles of estimation theory. It is also the reason why the
same algorithm circulates under different names in different communities: Tikhonov regularization
\cite{Tikhonov1977} and ridge regression \cite{Hoerl1970} in inverse problems and statistics, the
Wiener filter in signal processing \cite{Wiener1949,Kailath2000linear},
kriging in geostatistics \cite{Matheron1973,Myers1992}, and
Gaussian-process regression in machine learning \cite{Rasmussen2006,Bishop2006}. In finite
dimensions, the equivalence is elementary. Let us consider the linear measurement model
$\M y=\M H \M x + \V \epsilon$, where the signal $\M x \in \R^N$ is a realization of a Gaussian
random vector $\VRand X \sim \Spc N(\V 0,\M C_{\VRand X})$ and where
$\V \epsilon \sim \Spc N(\V 0,\sigma^2 \M I)$ is an independent measurement noise. The posterior
distribution of $\VRand X$ given $\M y$ is Gaussian, so that its mean and its mode coincide: the
minimum-mean-square-error (MMSE) estimator and the maximum-a-posteriori (MAP) estimator are one and
the same. The MAP estimator minimizes the negative log-posterior, which makes it the solution of a
penalized least-squares problem. What is specific to the Gaussian hypothesis is that
the penalty is quadratic and that the minimizer is also the MMSE estimator: the common estimator is
the solution of the Tikhonov problem
\begin{align}
\label{Eq:TikIntro}
\M x_{\lambda}=\arg \min_{\M x \in \R^N} \left( \|\M y - \M H \M x\|_2^2 + \lambda \|\M L \M x\|_2^2\right)
\end{align}
with $\lambda=\sigma^2$, provided that the regularization matrix $\M L$ whitens the prior, in the
sense that $\M L \M C_{\VRand X}\M L^\Tr=\M I$.
This classic correspondence is the template that the present paper transposes to the continuum. Our deterministic reference is the formulation of the variational theory of splines in \cite{Unser2026a,Unser2026b}, which is the most extensive that we are aware of: the regularization operator $\Lop$ is any spline-admissible operator, whose null space may be nontrivial, the signals are generalized functions, and the measurements are linear functionals that are continuous on the native space. Our aim is to establish the correspondence at this level of generality.
The generalized splines take the role of the linear estimators: they are the solutions of the
continuous-domain counterparts of \eqref{Eq:TikIntro}, and they are of the form
$f=\sum_{m=1}^M a_m \Op A\{\nu_m\}+p$, where the $\nu_m$ are the measurement functionals, $\Op A$ is
the Green's operator (conditional inverse) of the normal operator $\Lop^\ast\Lop$, and $p$ lies in
the null space of $\Lop$. The generalized Gaussian processes (GGPs) take the role of the Gaussian
random vectors; they are random
elements of the dual $\Spc S'$ of a nuclear Gelfand pair $(\Spc S,\Spc S')$, in the sense of Gelfand
and Vilenkin \cite{Gelfand-Vilenkin1964}; the prototype of such a pair is the Schwartz space of
rapidly decaying test functions together with its dual, the space of tempered distributions.

The continuous-domain version of this correspondence has a long history.
It was uncovered by Kimeldorf and Wahba, first for the case where $\Op L$ is shift-invariant with a trivial null space \cite{Kimeldorf1970} and then for the Tchebycheffian splines on an interval \cite{Kimeldorf1971}, and it was developed into a complete statistical methodology in \cite{Wahba1990,Silverman1985}. In parallel, Parzen
\cite{Parzen1961} recast time-series estimation in the language of reproducing-kernel Hilbert
spaces (RKHS) \cite{Aronszajn1950}, while Kailath and his students connected RKHS norms,
innovations representations, and Fredholm equations
\cite{Kailath1970,Kailath1972,Weinert1974,Weinert1978}. The modern kernel methods of machine
learning inherit the same duality \cite{Scholkopf2001,Hofmann2008,Kanagawa2018,Owhadi2019},
as does the stochastic-partial-differential-equation approach to spatial statistics
\cite{Lindgren2011,Lindgren2022}, in which a Gaussian field is specified through the operator
that whitens it rather than through its covariance.

There is, however, a loose end in this otherwise satisfying picture, and it affects precisely
the estimators that approximation theory favors. When the regularization operator $\Lop$ has a
nontrivial null space $\Spc N_{\V p}=\mathrm{span}\{p_1,\dots,p_{N_0}\}$---which is the case of
the derivative operators $\Dop^{m}$ that generate the polynomial splines, the workhorses of the
field---the $\Lop$-smoothing spline loses its MMSE optimality. What survives is the weaker
statement that it is the best {\em linear unbiased} estimator in a so-called fixed-effect model,
where the bias is measured relative to a fixed but unknown drift $p \in \Spc N_{\V p}$
\cite[Sec.\ 1.5]{Wahba1990}. Geostatistics meets the same obstacle: the intrinsic random functions
of order $k$ of Matheron \cite{Matheron1973} are defined only modulo the polynomials of degree $k$,
and universal kriging is, likewise, a best linear unbiased predictor. This is frustrating, for
these are the estimators with the best deterministic approximation rates, and they owe those
rates to the very null space that degrades their statistical status. The reason for the
degradation is easy to identify: the process whose covariance operator is the inverse of
$\Lop^\ast\Lop$ does not exist when $\Lop$ is not invertible, so that there is no prior with
respect to which a Bayesian estimator could be optimal.

The resolution that we propose is to construct the missing prior, in a setting that is broad enough to accommodate signals that are generalized functions. The conditioned parametrization of \cite{Unser2026b} provides the two ingredients: a right-inverse $\Op L^{-1}_{\V \phi}$ of $\Lop$ that is made unique by $N_0$ boundary conditions $\V \phi(f)=\big(\langle \phi_1,f\rangle,\dots,\langle \phi_{N_0},f\rangle\big)=\V 0$, and the associated Gram operator $\Op A_{\V \phi}=\Op L^{-1}_{\V \phi}\Op L^{-1\ast}_{\V \phi}$, which is symmetric and $\V p$-conditionally positive-definite: its quadratic form is positive on the test functions that annihilate $\Spc N_{\V p}$, while it vanishes on the boundary functionals $\phi_n$. The object that we are after is the Gaussian process
with covariance operator $\Op A_{\V \phi}: \Spc S \to \Spc S'$. Since $\Op A_{\V \phi}$ induces the norm of $\Spc H'_\Lop$ (Section \ref{Sec:NativeRecap}), it is not of trace class, and Sazonov's theorem \cite{Sazonov1958,Bogachev1998} rules out a realization of this process in the native space $\Spc H_\Lop$ itself.
The correct home is the dual of a nuclear space, where the
theorem of Minlos \cite{Minlos1959,Gelfand-Vilenkin1964} converts the continuity of a
characteristic functional into a bona fide probability measure. This is the same
$(\Spc S,\Spc S')$ sandwich that hosts the native space \cite{Unser2026a,Unser2026b}, which is
what makes the two constructions superposable.

Our contributions are the following.
\begin{enumerate}
\item {\em Generalized Gaussian processes on Hilbert--Gelfand spaces}
(Section \ref{Sec:ExtendedProcesses}).
Our starting point is the classical result of Gelfand and Vilenkin (Theorem \ref{Theo:GGaussProcess}), which attaches a GGP $\gRand g \sim \Spc N(\mu,\Op Q)$ to any symmetric positive-semi-definite covariance operator $\Op Q: \Spc S \to \Spc S'$. In Theorem \ref{Theo:GGPcondi}, we extend the action of the process to test functions in the dual $\Spc H'_\Lop \supset \Spc S$ of the native space, and we characterize the ``degenerate'' class that corresponds to a $\V p$-conditionally positive-definite covariance: such a process vanishes almost surely on the span $\Spc N_{\V \phi}$ of the boundary functionals, which is the probabilistic transcription of the boundary conditions.
\item {\em Innovation models} (Section \ref{Sec:GaussianSDE}).
We show that the processes of this class are exactly those that admit the innovation representation $\gRand G=\Op L^{-1}_{\V \phi}\{\gRand W\}$ in terms of a Gaussian white noise $\gRand W$, or equivalently, that solve the stochastic differential equation $\Lop \gRand G=\gRand W$ subject to $\V \phi(\gRand G)=\V 0$. The formalism subsumes the classical integral representations: the fractional Brownian motion of Mandelbrot and Van Ness \cite{Mandelbrot1968}
%, which was introduced in this very journal, 
is recovered as the case $\Lop=\Dop^{\gamma}$ with a single boundary condition, and its higher-order extensions are worked out in Appendix \ref{App:fBm}.
\item {\em Congruence} (Proposition \ref{Theo:Congruence}).
The Gaussian space of the process---the closed linear span of the random variables $\gRand G(\varphi)$---is isometric to the core $\Spc H_{\Op A}=\Lop^\ast(L_2)$, which is the part of the dual $\Spc H'_\Lop$ that does not depend on the boundary conditions. Under this isometry, the Riesz pairing of the native space is the covariance of the process, and the innovation provides an orthonormal coordinate system. This is the statement that welds the deterministic and the stochastic descriptions together; Table \ref{Tab:Dictionary} lists the correspondence term by term.
\item {\em Reproducing kernels and covariance functions}
(Section \ref{Sec:ReproKernel}).
For the classic setting where $\Spc H_\Op L$ is an RKHS on $\R^d$, we identify in Theorem \ref{Theo:ReproKernel} a single bivariate function $a_{\V \phi}$ as the Schwartz kernel of the Gram operator, the reproducing kernel of the primary native space $\Spc H_{\Lop,\V \phi}$ (the functions of $\Spc H_\Lop$ that satisfy the boundary conditions), and the covariance function of the process. The formula involves the Green's function of the normal operator $\Lop^\ast\Lop$, corrected in each variable by the finite-rank term that carries the boundary conditions. (For a trivial null space, the correction vanishes and $a_{\V \phi}$ reduces to the Green's function itself, which is the classical kernel \cite[Section~\ref{II-Sec:RKHS}]{Unser2026b}.) In Theorem \ref{Theo:OptoKernelPosDef}, we match the operator-level notion of ($\V p$-conditional) positive definiteness with the kernel-level one of the radial-basis-function literature \cite{Buhmann2003,Wendland2005}. As an illustration, we give in Section \ref{Sec:fBf} a complete and explicit characterization of the isotropic fractional Brownian fields of arbitrary order on $\R^d$---the innovation kernel, the stochastic integral representation \`a la Mandelbrot, and the pointwise covariance---which appears to be new beyond the classical range $H \in (0,1)$ of the Hurst exponent, the range of L\'evy's field and of its fractional extensions.

\item {\em Splines as MMSE estimators} (Section \ref{Sec:MMSESplines}).
For the measurement model $\M y=\V \nu(g)+\V \epsilon$, we identify in Theorem \ref{Theo:MMSE} the MMSE estimator of $\gRand G(\varphi)$ with the generalized smoothing spline that is constrained by the boundary conditions $\V \phi(g)=\V 0$ \cite{Unser2026b}, evaluated at $\varphi$, with $\lambda=\sigma^2$ and without any restriction to linear or to unbiased estimators. This closes the loose end: once the prior is the conditioned process, a generalized spline is the MMSE estimator outright, and not merely the best linear unbiased one. In Corollary \ref{Theo:Unser}, we give the noise-free counterpart: every generalized spline interpolant is the MMSE interpolator of a corresponding GGP.
\end{enumerate}

A single mathematical object underlies the whole construction, and it admits three readings. In
the language of functional analysis, it is a conditional Riesz pair $(\Op A,\Op B)$ over the
nuclear Gelfand pair $(\Spc S,\Spc S')$ \cite{Unser2026a}; in the language of approximation
theory, it is the native space $\Spc H_\Lop$ of the $\V p$-admissible operator $\Lop$, together
with its representer theorem \cite{Unser2026b}; and in the language of probability, which is the
one of the present paper, it is the covariance structure of a generalized Gaussian process. What
makes the third reading possible is the conditioned parametrization: the boundary conditions that
single out the right-inverse $\Op L^{-1}_{\V \phi}$ are exactly the ones that make the Gaussian
process well-defined.

The paper proceeds as follows. In Section \ref{Sec:Preliminaries}, we recall the
functional-analytic setting of \cite{Unser2026a,Unser2026b} and the Minlos--Bochner theorem. In
Section \ref{Sec:GGP}, we develop the theory of generalized Gaussian processes, from the
abstract definition to the innovation models and the congruence. In
Section \ref{Sec:ReproKernel}, we make the link with the classical theory: we identify the
reproducing kernel of the native space with the covariance function of the process, and we
relate the operator-level and kernel-level notions of conditional positive definiteness. In
Section \ref{Sec:MMSESplines}, we establish the equivalence between generalized splines and
MMSE estimators, first in the non-degenerate case, and then in full generality. The framework
is illustrated in Section \ref{Sec:fBf} with the characterization of the fractional Brownian
fields on $\R^d$, and in Appendix \ref{App:fBm} with the higher-order extensions of
fractional Brownian motion (fBm).

\section{Preliminaries}
\label{Sec:Preliminaries}

\subsection{Nuclear Gelfand Pairs and Hilbert--Gelfand Spaces}
\label{Sec:GelfandPairs}
We briefly recall the functional-analytic setting of \cite{Unser2026a,Unser2026b}, to which
we refer for the complete exposition.

\begin{definition}[Nuclear Gelfand pair]
\label{Def:NuclearGelfandPair}
The couple $(\Spc S,\Spc S')$ that is formed by a topological vector space $\Spc S$ and
its topological dual $\Spc S'$ is called a \emph{nuclear Gelfand pair} if
\begin{enumerate}
\item $\Spc S$ is nuclear and Montel
\cite{Grothendieck1955,Treves2006};
\item there exists a continuous and dense embedding $\iota:\Spc S\to\Spc S'$ that is
canonical, in the sense that
$\langle\iota\varphi,\psi\rangle=\langle\iota\psi,\varphi\rangle$ for all
$\varphi,\psi\in\Spc S$.
\end{enumerate}
Unless stated otherwise, $\Spc S'$ is endowed with its weak-$\ast$ topology.
\end{definition}
The prototypical example is the pair $\big(\Spc S(\R^d),\Spc S'(\R^d)\big)$ formed by the
Schwartz space of smooth and rapidly decaying test functions on $\R^d$ and its dual, the space
of tempered distributions. 

A Hilbert space $\Spc H$ is said to be a {\em Hilbert--Gelfand space}
if $\Spc S\embedD \Spc H\embedD\Spc S'$, where the symbol ``$\embedD$'' denotes a continuous and
dense embedding; such a space can be constructed by completion as
$\Spc H=\overline{(\Spc S,\|\cdot\|_{\Spc H})}$, and its dual $\Spc H'$ is a Hilbert--Gelfand
space as well, with the duality pairing that is compatible with that of $(\Spc S,\Spc S')$. The
baseline example is the pivot space $L_2=\overline{(\Spc S,\|\cdot\|_{L_2})}$.

\subsection{Native Spaces and Spline-Admissible Operators}
\label{Sec:NativeRecap}
We now recall, without proof, the operator-theoretic objects of
\cite{Unser2026a,Unser2026b} that the stochastic constructions of this paper rest upon. The regularization (or whitening) operator that guides the construction is $\Op L: \Spc S\to L_2$. By suitable completion, this operator induces the native space $\Spc H_\Op L$, which is such that $\Op L: \Spc H_\Op L \to L_2$ with finite dimensional null space
$\Spc N_\V p={\rm span}(p_1,\cdots,p_{N_0})=\{f \in \Spc H_\Lop: \Op L f=0\}$. The adjoint of $\Op L$ is
$\Op L^\ast: L_2 \to \Spc H_\Op A \subset \Spc H'_\Op L$ where $\Spc H'_\Op L$ is the topological dual of $\Spc H_\Op L$ and $\Spc H_\Op A:=\Op L^\ast(L_2)$.

Since the conditions under which this construction is valid are central to what follows, let
us spell them out explicitly. Let $\V p=(p_1,\dots,p_{N_0})\subset \Spc S'$ be a linearly
independent family with $\Spc N_{\V p}={\rm span}\{p_n\}_{n=1}^{N_0}$, and let
$\Spc S_{\V p}=\{\varphi \in \Spc S: \V p(\varphi)=\V 0\}$ be the space of test functions
that annihilate it. A continuous symmetric operator $\Op T: \Spc S \to \Spc S'$ induces the
energy seminorm $|\varphi|_{\Op T}=\langle \Op T\varphi,\varphi\rangle^{1/2}$, whose
{\em attainable null space}
\begin{align}
\label{Eq:AttainableNull}
\Spc N_{\Op T}=\big\{p \in \Spc S':\ p=\lim_{k \to \infty}\chi_k \mbox{ in } \Spc S'
\mbox{ for some } (\chi_k)\subset \Spc S \mbox{ with } |\chi_k|_{\Op T}\to 0\big\}
\end{align}
is the space of generalized functions that are reached as limits of test functions of
vanishing energy \cite[Definition~\ref{II-Def:AttainableNull}]{Unser2026b}. Let us stress
that $\Spc N_{\Op T}$ is generally non-trivial even when $|\cdot|_{\Op T}$ is a genuine norm
on $\Spc S$, since the limit is taken in the much weaker topology of $\Spc S'$.

\begin{definition}[$\V p$-admissible operator]
\label{Def:Admissible}
Let $\Lop: \Spc S \to L_2$ be a continuous linear operator that is injective and whose range
is dense in $L_2$, and let us assume that its normal operator
$\Op B=\Lop^\ast\Lop: \Spc S \to \Spc S'$ admits a {\em conditional inverse}
$\Op A: \Spc S \to \Spc S'$, in the sense that
\begin{align}
\label{Eq:CondInverse}
\Op A\Op B=\Identity \ \mbox{ on } \Spc S
\qquad \mbox{ and } \qquad
\Op B\Op A=\Identity \ \mbox{ on } \Spc S_{\V p}.
\end{align}
Then, $\Lop$ is said to be {\em $\V p$-admissible} if the family $\V p$ spans the attainable
null space of the spline energy $|\varphi|_{\Op B}=\|\Lop \varphi\|_{L_2}$; that is, if
\begin{align}
\label{Eq:NullMatching}
\Spc N_{\Op B}=\Spc N_{\V p}.
\end{align}
\end{definition}

The single condition \eqref{Eq:NullMatching} says that $\V p$ is neither too small nor too
large, and it is what makes the null space of the continuous extension
$\Lop: \Spc H_\Lop \to L_2$ be exactly $\Spc N_{\V p}$. Let us also point at the asymmetry
between the two operators, which is the delicate point of the construction. The normal
operator $\Op B$ is positive-definite on the whole of $\Spc S$, since
$|\varphi|_{\Op B}=\|\Lop \varphi\|_{L_2}$ and $\Lop$ is injective, while its conditional
inverse $\Op A$ is only {\em $\V p$-conditionally positive-definite}, in the sense that
$\langle \Op A\varphi,\varphi\rangle>0$ for all $\varphi \in \Spc S_{\V p}\backslash\{0\}$;
this is also why the second identity in \eqref{Eq:CondInverse} is restricted to
$\Spc S_{\V p}$. Under the hypotheses of Definition \ref{Def:Admissible}, the condition
\eqref{Eq:NullMatching} is equivalent to $(\Op A,\Op B)$ being a {\em $\V p$-conditional
Riesz pair} \cite[Corollary~\ref{I-Cor:GreenNormal}]{Unser2026a}: the two completions
$\Spc H_{\Op A}=\overline{(\Spc S_{\V p},|\cdot|_{\Op A})}$ and
$\dot{\Spc H}_{\Op B}=\overline{(\Spc S,|\cdot|_{\Op B})}$ are then put in duality by
$\Op A$ and $\Op B$, with $\Spc H_{\Op A}=\Lop^\ast(L_2)$ a genuine Hilbert subspace of
$\Spc S'$ and $\dot{\Spc H}_{\Op B}\cong L_2$ a space of equivalence classes modulo
$\Spc N_{\V p}$.

Let us now turn to the objects that the boundary conditions bring in. Given a $\V p$-admissible operator $\Lop$, there always exists a biorthogonal family $\V \phi=(\phi_1,\dots,\phi_{N_0}) \subset \Spc H_\Op L'$ such that $\langle p_m,\phi_n\rangle=\delta_{m,n}$ \cite[Lemma~\ref{I-Lem:Biorthogonal}]{Unser2026a}. To any such family
$\V \nu=(\nu_1,\dots,\nu_{N_0})$, we associate the rank-$N_0$ operator
$\Op R_{\V \nu}\eqdef\sum_{n=1}^{N_0} \nu_n\,\langle \nu_n,\cdot\rangle$. The relevant
projectors are then
\begin{align}
\label{III-Eq:Projectors}
\Op P_{\V \phi}=\Identity-\Op R_{\V \phi}\Op R_{\V p}: \Spc H'_\Op L \to \Spc H_\Op A
\end{align}
and its adjoint $\Op P^\ast_{\V \phi}=\Identity-\Op R_{\V p}\Op R_{\V \phi}: \Spc H_\Op L \to \Spc H_{\Lop,\V \phi}$.
The three objects that we shall use repeatedly are the following.
\begin{enumerate}
\item The {\em Gram operator}, which is symmetric and $\V p$-conditionally positive-definite, is defined as
\begin{align}
\label{III-Eq:Aphi}
\Op A_{\V \phi}=\Op P^\ast_{\V \phi}\Op A\Op P_{\V \phi}: \Spc S \to \Spc S'.
\end{align}
Its augmentation $\Op A_+=\Op A_{\V \phi}+\Op R_{\V p}$
is
positive-definite.
\item The {\em conditioned inverse} is the unique right-inverse of $\Lop$ that imposes the boundary conditions
$\V \phi(f)=\V 0$. It is given by
\begin{align}
\label{III-Eq:gInverse}
\Op L^{-1}_{\V \phi}=\Op A_{\V \phi}\Lop^\ast: L_2 \to \Spc S',
\end{align}
and it factorizes the Gram operator as \big(see \cite[Proposition~\ref{I-Prop:Factorization}]{Unser2026a}\big)
\begin{align}
\label{III-Eq:Factorization}
\Op A_{\V \phi}=\Op L^{-1}_{\V \phi}\Op L^{-1\ast}_{\V \phi}.
\end{align}
\item The {\em native space} $\Spc H_\Op L$ \cite[Theorem~\ref{II-Theo:NativeSpace}]{Unser2026b} is the Hilbert--Gelfand space  equipped with the inner product
\begin{align}
\label{III-Eq:InnerProd}
\langle f,g \rangle_{\Spc H_\Op L}=\langle \Op L f, \Op L g\rangle_{L_2}+\sum_{n=1}^{N_0}\langle \phi_n, f\rangle\langle \phi_n, g\rangle.
\end{align}
Consequently, $\Spc H_\Lop=\overline{(\Spc S,\|\cdot\|_{\Spc H_\Op L})}$ and it decomposes as
$\Spc H_\Lop=\Spc H_{\Lop, \V \phi} \oplus \Spc N_{\V p}$ with
$\Spc H_{\Lop, \V \phi}=\{f \in \Spc H_\Lop: \V \phi(f)=\V 0\}=\Op L^{-1}_{\V \phi}(L_2)$,
while its dual $\Spc H'_\Lop=\overline{(\Spc S,\|\cdot\|_{\Spc H'_\Lop})}$, with
$\|\varphi\|^2_{\Spc H'_\Lop}=\langle \Op A_+\varphi,\varphi\rangle$, decomposes as
$\Spc H'_\Lop=\Spc H_{\Op A} \oplus \Spc N_{\V \phi}$ with
$\Spc H_{\Op A}=\Lop^\ast(L_2)=\{\nu \in \Spc H'_\Lop: \V p(\nu)=\V 0\}$. The Riesz map
from $\Spc H_{\Op A}$ onto $\Spc H_{\Lop, \V \phi}$ is $\Op A_{\V \phi}$, so that
\begin{align}
\label{Eq:RieszPairingIII}
\langle \nu,\upsilon\rangle_{\Spc H'_\Lop}=\langle \Op A_{\V \phi} \nu,\upsilon\rangle
\qquad \mbox{ for all } \nu,\upsilon \in \Spc H_{\Op A}.
\end{align}
\end{enumerate}
We refer to $\Spc H_{\Op A}=\Lop^\ast(L_2)$ as the {\em core} of the construction: it is the
only member of the family that is independent of the boundary conditions $\V \phi$, from
which all the others are generated once $\V \phi$ is selected, and it is the space that
carries the second-order structure of the process (see Proposition \ref{Theo:Congruence}, which shows that $\Spc H_\Op A$ is congruent to the Gaussian space
$\Spc G_{\gRand G}$).

The theme of the present paper is that the Gram operator $\Op A_{\V \phi}$ plays a double role. In the variational theory of \cite{Unser2026a,Unser2026b}, it is the Riesz map of the core $\Spc H_{\Op A}$, which induces the inner product \eqref{Eq:RieszPairingIII}. Here, it is the covariance operator of a generalized Gaussian process $\gRand G \sim \Spc N(0,\Op A_{\V \phi})$ (Theorem \ref{Theo:GGaussProcess}), so that \eqref{Eq:RieszPairingIII} becomes the identity $\langle \nu,\upsilon\rangle_{\Spc H'_\Lop}={\rm Cov}_{\gRand G}(\nu,\upsilon)$ and every term of the variational construction acquires a stochastic counterpart. Table \ref{Tab:Dictionary} summarizes the correspondence, whose
mathematical content is the congruence of Proposition \ref{Theo:Congruence}.

\begin{table}[t]
\caption{Dictionary between the deterministic objects of the variational theory
\cite{Unser2026a,Unser2026b} and their stochastic counterparts in the present paper.}
\label{Tab:Dictionary}
\centering
\vspace{1ex}\setlength{\tabcolsep}{4pt}
\begin{tabular}{lll}
\hline\\[-2ex]
Variational theory & This paper & Role \\[1ex]
\hline\\[0.ex]
$\Lop$ ($\V p$-admissible) & $\Lop$ & whitening operator \\[1ex]
$\Op A=(\Op L^\ast\Op L)^{-1}$ & $ \Op A$ & Conditional inverse of $(\Op L^\ast\Op L)$ \\[1ex]
$\Op A_{\V \phi}=\Op P^\ast_{\V \phi}\Op A\Op P_{\V \phi}$ & $\Op C_{\Lop,\V \phi}$ & covariance operator \\[1ex]
$\Op A_+=\Op A_{\V \phi}+\Op R_{\V p}$ & $\Op A_+$ & augmented (proper) covariance \\[1ex]
$\Op L^{-1}_{\V \phi}=\Op A_{\V \phi}\Lop^\ast$ & shaping operator $\Op T$ & synthesis of $\gRand G$ from $\gRand W$ \\[1ex]
$\Op A_{\V \phi}=\Op L^{-1}_{\V \phi}\Op L^{-1\ast}_{\V \phi}$ & $\Op C_{\Lop,\V \phi}=\Op T\Op T^\ast$ & innovation model \\[1ex]
$\V \phi(f)=\V 0$ & $\V \phi(\gRand G)=\V 0$ & boundary conditions \\[1ex]
$\Spc H_\Lop$ (native space) & home of $g_\lambda$ & MMSE estimate \\[1ex]
$\Spc H'_\Lop$ (dual) & home of $\varphi$ and $\nu_m$ & test/measurement functionals \\[1ex]
$\Spc H_{\Op A}=\Lop^\ast(L_2)$ & $\Spc H_{\Op A}$ & core, congruent to $\Spc G_{\gRand G}$ \\[1ex]
$\Spc H_{\Lop,\V \phi}$ & $\Spc H_{\Lop,\V \phi}$ & primary native space \\[1ex]
$\langle \nu,\upsilon\rangle_{\Spc H'_\Lop}$ & ${\rm Cov}\big(\gRand G(\nu),\gRand G(\upsilon)\big)$ & Riesz pairing $=$ covariance \\[1ex]
$a_{\V \phi}(\V x,\V y)$ & $\Exp\{\gRand G(\V x)\gRand G(\V y)\}$ & reproducing kernel $=$ covariance \\[1ex]
$L_2$ (via $\Lop$) & innovation $\gRand W$ & orthonormal coordinates \\[1ex]
\hline
\end{tabular}
\end{table}

\subsection{Probability Measures on $\Spc S'$: The Minlos Theorem}
\label{Sec:Minlos}
The stochastic part of our story requires probability measures on infinite-dimensional
topological vector spaces. The adequate device to specify them is the {\em characteristic
functional} $\ProbF_\gRand g: \Spc S \to \C$ with
$\ProbF_\gRand g(\varphi)=\Exp\{\ee^{\jj \langle \gRand g, \varphi\rangle}\}$, which is the
infinite-dimensional generalization of the characteristic function
\eqref{Eq:MultiGaussFourier} of a random vector. In finite dimensions, Bochner's theorem
asserts that the characteristic functions on $\R^N$ are exactly the continuous
positive-definite functions that are normalized at the origin, and that they are in one-to-one
correspondence with the probability measures on $\R^N$. The fundamental result of Minlos
\cite{Minlos1959} extends this correspondence to any nuclear space, with the measure being
supported by the dual \cite{Gelfand-Vilenkin1964}.

\begin{theorem}[Minlos--Bochner]
\label{Theo:Minlos}
Let $\Spc S$ be a real nuclear locally convex space and $\Spc S'$ its topological dual
(with the cylindrical $\sigma$-algebra).
Let $\ProbF : \Spc S \to \C$ be a map that satisfies the three defining properties of a {\em characteristic functional}:

\begin{enumerate}
\item Normalization: $\ProbF(0)=1$,
\item Positive-definiteness: For all $N\in\mathbb N$,
      $\varphi_1,\dots,\varphi_N\in \Spc S$ and $c_1,\dots,c_N\in\C$,
\begin{align}
\label{Eq:PSfunction}
\sum_{m=1}^N \sum_{n=1}^N c_m \overline{c_n}\,
\ProbF (\varphi_m-\varphi_n) \ge 0,
\end{align}
\item Continuity: $\ProbF$ is continuous at $0$.
\end{enumerate}

Then, there exists a unique probability measure $\Prob$ on $\Spc S'$ such that
\begin{align}
\ProbF(\varphi)
=\int_{\Spc S'} \ee^{\jj\langle \omega,\varphi\rangle}\, \dint \Prob(\omega),
\qquad \varphi\in \Spc S.
\end{align}

Conversely, every probability measure on $\Spc S'$ has a characteristic
functional that satisfies {\rm (1)--(3)}.

\end{theorem}
An important complement is that the three defining properties in Theorem \ref{Theo:Minlos} imply that $\ProbF$ is continuous on $\Spc S$, Hermitian symmetric (i.e., $\ProbF(\varphi)=\overline{\ProbF(-\varphi)}$), and such that $0\le |\ProbF(\varphi)|\le \ProbF(0)=1$ for any $\varphi \in \Spc S$. The first of these implications explains why
the third defining hypothesis is often replaced by the more natural requirement that $\ProbF$ be continuous on $\Spc S$.

Let us explain why the dual of a nuclear space---rather than a Hilbert space---is the correct
home for such measures. The naive transposition of Bochner's theorem to an infinite-dimensional
Hilbert space fails, and it fails on the most important example: the functional
\begin{align}
\label{Eq:WhiteNoiseCF}
\ProbF_{\gRand w}(\varphi)=\ee^{-\frac{1}{2}\|\varphi\|^2_{L_2}},
\end{align}
which is continuous and positive-definite on $L_2$, and which is the natural candidate for the
characteristic functional of a white Gaussian noise with correlation
$\Exp\{\langle \gRand w,\varphi\rangle \langle \gRand w,\psi\rangle\}=\langle \varphi,\psi\rangle$.
Were $\gRand w$ a random element of $L_2$, its coefficients
$\wRand_k=\langle \gRand w, e_k\rangle$ in an orthonormal basis $(e_k)_{k \in \N}$ would be
i.i.d.\ $\Spc N(0,1)$, so that $\|\gRand w\|^2_{L_2}=\sum_{k} |\wRand_k|^2=+\infty$ almost
surely: the candidate measure assigns zero mass to the very space that is supposed to support
it.
%, and what survives on $L_2$ is only a finitely additive cylinder measure. 
Sazonov's theorem
\cite{Sazonov1958,Bogachev1998} pinpoints the obstruction: on a Hilbert space, a Gaussian
functional $\ee^{-\frac{1}{2}\langle \Op Q \varphi,\varphi\rangle}$ is the characteristic
functional of a countably additive measure if and only if its covariance $\Op Q$ is
trace-class, which the identity operator in \eqref{Eq:WhiteNoiseCF} is not in infinite
dimensions.

The theorem of Minlos resolves the difficulty by enlarging the target space. When $\Spc S$ is
nuclear, the canonical maps between the Hilbertian seminorm completions that define its
topology are themselves nuclear (trace-class), which means that the correction demanded by
Sazonov is built into the topology of $\Spc S$. Consequently, the mere continuity of $\ProbF$ on
$\Spc S$ suffices to produce a probability measure, albeit on the larger space $\Spc S'$. In
particular, white Gaussian noise exists as a generalized process with realizations in
$\Spc S'$, but $\Prob_{\gRand w}(L_2)=0$. This is the division of labor within the sandwich
$\Spc S \embedD \Spc H_\Lop \embedD \Spc S'$: the native Hilbert space hosts the
deterministic reconstruction, while $\Spc S'$ hosts the stochastic model, whose realizations
almost surely fall outside the native space $\Spc H_\Lop$ (and, for white noise, outside $L_2$).

\section{Generalized Gaussian Processes}
\label{Sec:GGP}

\subsection{Abstract Generalized Gaussian Processes}
\label{Sec:GGP0}
The definition of abstract Gaussian processes is motivated by two important observations: the first is that any linear functional of a Gaussian
random vector is a scalar Gaussian random variable. The second is that the characteristic function of the Gaussian random vector $\VRand X\sim\Spc N(\V \mu_{\VRand X},\M C_{\VRand X})$ with mean vector $\V \mu_{\VRand X}\in \R^N$ and covariance matrix  $\M C_{\VRand X}\in \R^{N \times N}$ is given 
by 
\begin{align}
\label{Eq:MultiGaussFourier}
\hat p_{\VRand X}(\V \xi)= \Exp\{\ee^{\jj \langle \V \xi, \VRand X\rangle}\}=\exp\left( \jj \langle \V \mu_{\VRand X},\V \xi\rangle - \tfrac{1}{2} \langle \M C_{\VRand X}\V \xi,\V \xi\rangle \right).
\end{align}
This function is the conjugate Fourier transform of the corresponding multivariate Gaussian probability density function (pdf). 

Generalized Gaussian processes take one more step of abstraction than the classical random
functions: they are defined either as random elements of the dual $\Spc S'$ of a nuclear
space $\Spc S$, or as random linear functionals acting on $\Spc S$. The choice
$\Spc S=\Spc S(\R^d)$ is the one of interest here, but the framework covers every Gaussian
object of interest, down to random vectors ($\Spc S=\R^N$) and discrete signals
($\Spc S=\Spc S(\Z)$).

To explain the concept, we recall that a generalized function $g \in \Spc S'$ is nothing
but a rule that associates the real number $g(\varphi)=\langle g,\varphi\rangle$ with every
test function $\varphi \in \Spc S$, linearly and continuously.
A generalized stochastic process randomizes this description \cite{Gelfand-Vilenkin1964,Bierme2018,Dalang2017}, which can be done from two complementary points of view: one either randomizes the {\em object} $g$, or the {\em rule} itself.
\begin{definition}[Generalized stochastic process]
\label{Theo:GGauss}
A {\em generalized stochastic process} in $\Spc S'$ admits the two following equivalent
descriptions.
\begin{enumerate}
\item {\em Random element of $\Spc S'$} (randomization of the object): $\gRand g$ is an
$\Spc S'$-valued random variable $\gRand g: \Omega \to \Spc S'$ that is defined on some
probability space $(\Omega,\Spc F,\Prob)$; equivalently, it is specified by a probability
measure $\Prob_{\gRand g}$ on $\Spc S'$, according to which the process produces its
{\em realizations} $g=\gRand g(\omega) \in \Spc S'$, $\omega \in \Omega$.
\item {\em Random linear functional on $\Spc S$} (randomization of the rule):
$\gRand g=\{\gRand g(\varphi)\}_{\varphi \in \Spc S}$ is a collection of real-valued
random variables indexed by $\Spc S$ such that
\begin{itemize}
\item[(i)] (linearity)
$\gRand g(a_1 \varphi_1 + a_2 \varphi_2) = a_1 \gRand g(\varphi_1) + a_2 \gRand g(\varphi_2)$
almost surely, for any $\varphi_1,\varphi_2 \in \Spc S$ and $a_1,a_2 \in \R$;
\item[(ii)] (consistency) the joint laws of the vectors
$\big(\gRand g(\varphi_1),\allowbreak\dots,\allowbreak\gRand g(\varphi_N)\big)$, with
$\varphi_1,\dots,\varphi_N \in \Spc S$ and $N \in \N$, fit together coherently, in the
sense that they agree with each other under reordering and under restriction to subsets;
\item[(iii)] (continuity) $\varphi_k \to \varphi$ in $\Spc S$ implies that
$\gRand g(\varphi_k) \to \gRand g(\varphi)$ in probability.
\end{itemize}
\end{enumerate}
The process is {\em Gaussian}---a generalized Gaussian process (GGP)---if the random
vector $\VRand Y=\big(\gRand g(\varphi_1), \dots, \gRand g(\varphi_N)\big)$ has a
multivariate Gaussian distribution for every finite collection
$\varphi_1,\dots,\varphi_N \in \Spc S$.
\end{definition}

The direction $1 \Rightarrow 2$ of the equivalence is elementary: the variables
$\gRand g(\varphi)=\langle \gRand g,\varphi\rangle$ inherit the linearity and the
continuity from the realizations pointwise (hence almost surely and in probability), and
their joint laws are automatically coherent, because they are all marginals of the common
measure $\Prob_{\gRand g}$. The converse is the fundamental benefit of the nuclear
setting. By the Kolmogorov extension theorem, the coherent family of Description 2
specifies a cylindrical measure, whose characteristic functional
$\varphi \mapsto \Exp\{\ee^{\jj \gRand g(\varphi)}\}$ is normalized, positive-definite
\big(by Item (i), as in \eqref{Eq:PDzn}\big), and continuous on $\Spc S$ \big(by
Item (iii) and dominated convergence\big); the Minlos--Bochner theorem
(Theorem \ref{Theo:Minlos}) then upgrades that cylindrical measure to a genuine probability
measure $\Prob_{\gRand g}$ on $\Spc S'$ with the same finite-dimensional laws
\cite{Gelfand-Vilenkin1964}. In that case, one may always take the canonical probability
space $(\Omega,\Spc F,\Prob)=(\Spc S',\Spc B,\Prob_{\gRand g})$, where $\Spc B$ is the
cylindrical $\sigma$-field, with the tautological action
$\gRand g(\varphi): \omega \mapsto \langle \omega,\varphi\rangle$. Description 2 is the operative one in this paper: it is the
form in which we construct the processes in Theorem \ref{Theo:GGaussProcess}, and the one
whose domain we extend beyond $\Spc S$ in Section \ref{Sec:ExtendedProcesses}.
Description 1 is what makes statements about the realizations meaningful, including the characterization of their MMSE estimators in Section \ref{Sec:MMSESplines}.

Since Gaussian random variables have a finite variance, all the observations
$\gRand g(\varphi)$ of a GGP live in the Hilbert space $L_2(\Omega,\Prob)$ of (equivalence
classes of) random variables with finite second moments over the probability space
$(\Omega,\Spc F,\Prob)$ of Definition \ref{Theo:GGauss}, which is equipped with the inner
product
$\langle \Rand X,\Rand Y\rangle_{L_2(\Omega,\Prob)}=\Exp\{\Rand X\hspace{0.08em}\Rand Y\}$
and in which convergence is convergence in the quadratic mean. This space---not to be
confused with $L_2=L_2(\R^d)$---is the ambient structure for the second-order theory of the
paper: the extension of Section \ref{Sec:ExtendedProcesses} is obtained by
$L_2(\Omega,\Prob)$-continuity, and the Gaussian space of
Proposition \ref{Theo:Congruence} is one of its closed subspaces.

As in the finite-dimensional setting, a GGP is characterized by its first and second-order moments, with the difference that the mean $\mu=\Exp\{ \gRand g\}\in \Spc S'$
is now a continuous linear functional on $\Spc S$, while the ``covariance matrix'' takes the form of the symmetric positive-semi-definite operator $\Op Q: \Spc S \to \Spc S'$ \cite{Gelfand-Vilenkin1964}. These enable the determination
of the means and covariances of $\gRand g$ according to the rules
\begin{align*}
\mu(\varphi)&=\Exp\{\gRand g(\varphi)\}=\langle \mu,\varphi\rangle \\
{\rm Cov}_{\gRand G}(\varphi_1,\varphi_2\big)&=\Exp\Big\{\big(\gRand g(\varphi_1)-\mu(\varphi_1)\big)\big(\gRand g(\varphi_2)-\mu(\varphi_2)\big)\Big\}=\langle \Op Q \varphi_1,\varphi_2 \rangle.
\end{align*}
The corresponding generalized Gaussian process is denoted $\gRand g \sim \Spc N(\mu,\Op Q)$, in direct analogy with the finite-dimensional convention.

In the next theorem, which is due to Gelfand and Vilenkin \cite{Gelfand-Vilenkin1964}, we
recall that such a generalized Gaussian process exists for any $\mu \in \Spc S'$ and
symmetric, positive semi-definite operator $\Op Q$, and that it is characterized explicitly
by its {\em characteristic functional} on $\Spc S$---the infinite-dimensional counterpart
of \eqref{Eq:MultiGaussFourier}. We state the result in the form that the rest of the paper
relies upon, and we include its short proof, on which Theorem \ref{Theo:GGPcondi} builds.
In particular, the generalized Gaussian process $\gRand W \sim \Spc N(0,\Identity)$ 
is called a {\em white noise} or {\em Gaussian innovation}.

\begin{theorem}[Generalized Gaussian process in $\Spc S'$ {\cite{Gelfand-Vilenkin1964}}]
\label{Theo:GGaussProcess}
For any given $\mu \in \Spc S'$ (mean) and
symmetric, positive semi-definite operator $\Op Q: \Spc S \to \Spc S'$ (covariance), there exists a corresponding
generalized Gaussian process $\gRand g \sim \Spc N(\mu,\Op Q)$ on $\Spc S$ whose characteristic functional
is \begin{align}
\label{Eq:CharacGauss}
\ProbF_{\gRand g}(\varphi)=
\Exp\Big\{\ee^{ \jj \gRand g(\varphi)}\Big\} = \exp\left( \jj \langle \mu,\varphi\rangle - \tfrac{1}{2} \langle \Op Q \varphi, \varphi\rangle \right).
\end{align}
The latter specifies a Gaussian probability measure $\Prob_{\gRand g}$ on $\Spc S'$
that induces compatible families of multivariate Gaussian distributions: For each finite set $\{\varphi_1, \dots, \varphi_N\} \subset \Spc S$, we have that $\VRand Y={\big(\gRand g(\varphi_1), \dots,\gRand g(\varphi_N)\big)} \sim \Spc N(\V \mu_{\VRand Y},\M C_{\VRand Y})$ whose mean vector $\V \mu_{\VRand Y}\in \R^N$ and covariance matrix $\M C_{\VRand Y} \in \R^{N \times N}$ are given by
\begin{align}
[\V \mu_{\VRand Y}]_n=\langle \mu,\varphi_n\rangle
\label{Eq:Mean}
\\
[\M C_{\VRand Y}]_{m,n}={\rm Cov}_{\gRand g}(\varphi_m,\varphi_n)=\langle \Op Q \varphi_m,\varphi_n \rangle
\label{Eq:Cov}.
\end{align}
\end{theorem}

\begin{proof}
By the Minlos--Bochner theorem (Theorem \ref{Theo:Minlos}), it suffices to verify that the
right-hand side of \eqref{Eq:CharacGauss} is normalized, positive-definite \big(in the sense
of \eqref{Eq:PSfunction}\big), and continuous on $\Spc S$. The normalization holds by
inspection. The continuity holds because $\varphi_k \to \varphi$ in $\Spc S$ implies that
$\langle \mu,\varphi_k\rangle \to \langle \mu,\varphi\rangle$ and that
$\Op Q\varphi_k \to \Op Q\varphi$ in $\Spc S'$, while the duality pairing is hypocontinuous,
so that $\langle \Op Q\varphi_k,\varphi_k\rangle \to \langle \Op Q\varphi,\varphi\rangle$.
As for the positive-definiteness, it only involves finitely many test functions at a time.
Given $\varphi_1,\dots,\varphi_N \in \Spc S$, the vector $\V \mu_{\VRand Y}$ and the matrix
$\M C_{\VRand Y}$ of \eqref{Eq:Mean}--\eqref{Eq:Cov} specify a Gaussian vector
$\VRand Y \sim \Spc N(\V \mu_{\VRand Y},\M C_{\VRand Y})$, which exists because
$\M C_{\VRand Y}$ is symmetric and positive semi-definite (degenerate if it is singular). By
the bilinearity of the duality product,
$\langle \mu,\varphi_m-\varphi_n\rangle=[\V \mu_{\VRand Y}]_m-[\V \mu_{\VRand Y}]_n$ and
$\langle \Op Q(\varphi_m-\varphi_n),\varphi_m-\varphi_n\rangle
=(\V e_m-\V e_n)^\Tr \M C_{\VRand Y}(\V e_m-\V e_n)$, so that
$\ProbF_{\gRand g}(\varphi_m-\varphi_n)=\Exp\{\ee^{\jj(\Rand Y_m-\Rand Y_n)}\}$ by
\eqref{Eq:MultiGaussFourier}. Consequently, for all $c_1,\dots,c_N \in \C$,
\begin{align}
\label{Eq:PDzn}
\sum_{m=1}^N \sum_{n=1}^N c_m \overline{c_n}\, \ProbF_{\gRand g}(\varphi_m-\varphi_n)
=\Exp\Big\{\Big|\sum_{m=1}^N c_m\, \ee^{\jj \Rand Y_m}\Big|^2\Big\} \ge 0,
\end{align}
which proves the existence of the measure $\Prob_{\gRand g}$. Finally, the characteristic
function of $\big(\gRand g(\varphi_1),\dots,\gRand g(\varphi_N)\big)$ at
$\bwRand \in \R^N$ is, by the linearity of $\gRand g$,
$\Exp\big\{\ee^{\jj \gRand g(\sum_n \wRand_n\varphi_n)}\big\}
=\ProbF_{\gRand g}\big(\textstyle\sum_n \wRand_n \varphi_n\big)
=\exp\big(\jj \V \mu_{\VRand Y}^\Tr \bwRand-\tfrac{1}{2}\bwRand^\Tr \M C_{\VRand Y}\bwRand\big)$,
which is \eqref{Eq:MultiGaussFourier}; that is,
$\VRand Y \sim \Spc N(\V \mu_{\VRand Y},\M C_{\VRand Y})$.
\end{proof}
Theorem \ref{Theo:GGaussProcess} only requires the covariance operator to be positive
semi-definite. We call the GGP {\em proper} when $\Op Q$ is (strictly) positive-definite, and
{\em degenerate} (or improper) otherwise, in which case the measure $\Prob_{\gRand g}$ is
supported on a proper subspace of $\Spc S'$.

\subsection{Gaussian Processes Acting on Hilbert--Gelfand Spaces}
\label{Sec:ExtendedProcesses}
We now investigate the possibility of extending the action of a generalized Gaussian process in $\Spc S'$ to test functions $\varphi$ that may lie outside $\Spc S$.
This will be achieved by establishing the connection with the Hilbert--Gelfand spaces of \cite{Unser2026a,Unser2026b}, as recalled in Section \ref{Sec:NativeRecap}. The relevant class of covariance operators is that of the Gram operators $\Op A_{\V \phi}$ of \eqref{III-Eq:Aphi}, whose augmentation $\Op A_+=\Op A_{\V \phi}+\Op R_{\V p}$ underlies the dual pair $(\Spc H'_\Lop,\Spc H_\Lop)$, with $\Spc H'_\Lop=\Spc H_{\Op A} \oplus \Spc N_{\V \phi}$.
Note that this setting also covers the simpler non-degenerate case with $\V p=\{0\}$, which corresponds to the dual pair of Hilbert--Gelfand spaces $(\Spc H'_\Lop,\Spc H_\Lop)$ defined in \cite[Theorem~\ref{I-Theo:HilbertPositiveDef}]{Unser2026a}.

\begin{theorem}[GGP linked to a Hilbert--Gelfand space]
\label{Theo:GGPcondi}
Let $\Op A: \Spc S \to \Spc S'$ be a continuous linear operator that is $\V p$-conditionally positive-definite with corresponding $\Op A_\V \phi: \Spc H'_\Lop \to \Spc H_\Lop$ as recalled in Section \ref{Sec:NativeRecap}.
Then, the characteristic functional
\begin{align}
\label{Eq:CharacGauss2}
\varphi \mapsto \ProbF_{\gRand g}(\varphi)
= 
\exp\left( - \tfrac{1}{2} \langle \Op A_\V \phi \varphi, \varphi\rangle \right),\end{align}
specifies a generalized Gaussian process  $\gRand g \sim \Spc N(0,\Op A_{\V \phi})$ in $\Spc S'$ whose action is extendable to test functions $\varphi \in \Spc H'_\Lop\supset \Spc S$ where $\Spc H'_\Lop=\Spc H_{\Op A} \oplus \Spc N_{\V \phi}$
is the Hilbert space induced by the augmented operator $\Op A_+=\Op A_\V \phi+ \Op R_\V p$. Specifically, for each finite set $\{\varphi_1, \dots, \varphi_N\} \subset \Spc H'_\Lop$, 
this process induces compatible families of random Gaussian vectors $\VRand Y={\big(\gRand g(\varphi_1), \dots,\gRand g(\varphi_N)\big)} \sim \Spc N(\M 0,\M C_{\VRand Y})$ whose covariance matrices $\M C_{\VRand Y} \in \R^{N \times N}$ are given by
\begin{align}
[\M C_{\VRand Y}]_{m,n}=\langle \Op A_\V \phi \varphi_m,\varphi_n \rangle=\langle \Op A_\V \phi \widetilde{\varphi}_m,\widetilde{\varphi}_n \rangle=\langle \Op A \widetilde{\varphi}_m,\widetilde{\varphi}_n \rangle
\label{Eq:Cov2}
\end{align}
where $\widetilde{\varphi}_n=\Op P_{\V \phi}\{\varphi_n\} \in \Spc H_{\Op A}$ with $\Op P_{\V \phi}=\Identity- \Op R_\V \phi \Op R_\V p$. 
Moreover, $\VRand Y=\widetilde{\VRand Y}={\big(\gRand g(\widetilde{\varphi}_1), \dots,\gRand g(\widetilde{\varphi}_N)\big)}$ {\em almost surely}, which also indicates that $\gRand g$ vanishes on $\Spc N_\V \phi$ (degeneracy). 

Finally, the matrices $\M C_{\VRand Y}=\M C_{\widetilde{\VRand Y}}$ are positive-definite {\em if and only if} the corresponding $\widetilde{\varphi}_1,\dots,\widetilde{\varphi}_N \in \Spc H_{\Op A}$ are linearly independent.
For the non-degenerate scenario $\V p=\{0\}$, the latter property holds for any combination of linearly independent $\varphi_1,\dots,\varphi_N \in \Spc H'_\Lop \supset \Spc S$.

\end{theorem}

\begin{proof}
\medskip
\emph{Step 1: existence}: The Gram operator $\Op A_{\V \phi}: \Spc S \to \Spc S'$ is
continuous and symmetric, and it is positive semi-definite because
$\langle \Op A_{\V \phi}\varphi,\varphi\rangle=\|\Op L^{-1\ast}_{\V \phi}\varphi\|^2_{L_2}\ge 0$
by \eqref{III-Eq:Factorization}. Theorem \ref{Theo:GGaussProcess} with $\mu=0$ and
$\Op Q=\Op A_{\V \phi}$ therefore provides the GGP $\gRand g \sim \Spc N(0,\Op A_{\V \phi})$
in $\Spc S'$ with characteristic functional \eqref{Eq:CharacGauss2}, whose finite-dimensional
distributions over $\Spc S$ are the zero-mean Gaussians with covariance
$[\M C_{\VRand Y}]_{m,n}=\langle \Op A_{\V \phi}\varphi_m,\varphi_n\rangle$, which is the
first expression in \eqref{Eq:Cov2}.

\medskip
\emph{Step 2: the algebra of the conditioned covariance}: The key property here is that  $(\varphi,\psi) \mapsto \langle \Op A_\V \phi \varphi, \psi\rangle$ specifies a semi-inner product on $\Spc H'_\Lop$ (resp., an inner product on $\Spc H_\Op A$). Moreover, for all $\widetilde{\varphi},\widetilde{\psi} \in \Spc H_{\Op A}=\{\nu \in \Spc H'_\Lop: \V p(\nu)=\V 0\}$, we have that
\begin{align}
\label{Eq:EqinnerProdHA}
\langle \Op A_{\V \phi} \widetilde{\varphi}, \widetilde{\psi}\rangle=  \langle \Op A \widetilde{\varphi}, \widetilde{\psi}\rangle \mbox{ for any } \V \phi \subset \Spc H'_\Lop \mbox{ s.t. } \langle p_m,\phi_n\rangle=\delta_{m,n},
\end{align}
which is a direct consequence of the definition \eqref{III-Eq:Aphi} of the Gram operator, together with the property that $\Op P_{\V \phi}$ acts as the identity on $\Spc H_{\Op A}=\{\nu \in \Spc H'_\Lop: \V p(\nu)=\V 0\}$ \big(see \eqref{III-Eq:Projectors}\big).

\medskip
\emph{Step 3: extension to $\Spc H'_\Lop$}: The extension is obtained by continuity in the quadratic mean. For $\varphi \in \Spc S$,
the variance of the Gaussian random variable $\gRand g(\varphi)$ is
$\Exp\{|\gRand g(\varphi)|^2\}=\langle \Op A_{\V \phi}\varphi,\varphi\rangle
\le \langle \Op A_+\varphi,\varphi\rangle=\|\varphi\|^2_{\Spc H'_\Lop}$, where the inequality
holds because $\langle \Op R_{\V p}\varphi,\varphi\rangle=|\V p(\varphi)|^2_2\ge 0$. The map $\varphi \mapsto \gRand g(\varphi)$ is therefore a linear contraction from $\big(\Spc S,\|\cdot\|_{\Spc H'_\Lop}\big)$ into $L_2(\Omega,\Prob)$; that is, a linear map that does not increase the norm, with equality when $\V p(\varphi)=\V 0$. Since $\Spc H'_\Lop$ is
the completion of $\Spc S$ for that norm (Section \ref{Sec:NativeRecap}),
the map extends uniquely to a linear contraction $\Spc H'_\Lop \to L_2(\Omega,\Prob)$, and the
covariance identity \eqref{Eq:Cov2} passes to the limit, with its right-hand side given by
\eqref{Eq:EqinnerProdHA}. The Gaussianity of the finite-dimensional distributions is
preserved as well, because a limit in the quadratic mean of jointly Gaussian vectors is
jointly Gaussian.

\medskip
\emph{Step 4: degeneracy and non-degeneracy}: The direct-sum structure of $\Spc H'_\Lop$ implies that every $\varphi \in \Spc H'_\Lop$ has a unique decomposition as $\varphi=\widetilde{\varphi} + \phi$ with
$\widetilde{\varphi}=\Op P_{\V \phi}\{\varphi\} \in \Spc H_{\Op A}$ and $\phi=\varphi-\widetilde{\varphi}
\in \Spc N_{\V \phi}$. By linearity of the GGP, this yields $\Rand Y=\gRand g(\varphi)=\gRand g(\widetilde \varphi)+\gRand g(\phi)$
where $\gRand g(\widetilde \varphi)=\widetilde{\Rand Y}\sim \Spc N(0,\sigma^2_\Rand Y)$ with
$\sigma^2_{\Rand Y}=\langle \Op A_{\V \phi} \varphi,\varphi\rangle$ \big(from \eqref{Eq:Cov2}\big) and $\gRand g(\phi)=\Spc N(0,0)$ since $\langle \Op A_{\V \phi} \phi,\phi\rangle=0$
for all $\phi\in \Spc N_\V \phi$. The probabilistic interpretation of these properties is $\Rand Y=\widetilde{\Rand Y}$ and $\gRand g(\phi)=0$ (almost surely) for all $(\varphi,\phi) \in \Spc H'_\Lop \times \Spc N_\V \phi$. The final statement about the positive-definiteness of $\M C_{\VRand Y}$ follows from the positive-definiteness of $\Op A_\V \phi$ on $\Spc H_\Op A$, since $\V z^\Tr \M C_{\VRand Y}\V z=\big\langle \Op A_{\V \phi}\big\{\textstyle\sum_n z_n \widetilde{\varphi}_n\big\},\textstyle\sum_n z_n \widetilde{\varphi}_n\big\rangle$ vanishes exactly when $\sum_n z_n \widetilde{\varphi}_n=0$.

\end{proof}
The functional connection with Hilbert spaces in Theorem \ref{Theo:GGPcondi} should not be
over-interpreted: the fact that $\gRand g(\varphi)$ is a well-defined random variable for
every $\varphi \in \Spc H'_\Lop$ does not make $\gRand g$ a random process in the predual
$\Spc H_\Lop$, nor does it turn $\Prob_{\gRand g}$ into a measure on $\Spc H_\Lop$---the
white noise of Section \ref{Sec:Minlos} is again the witness. The rigorous home of these
processes remains the nuclear dual pair $(\Spc S,\Spc S')$, even when Hilbert spaces sit in
the background.

An attractive aspect of the extended framework of Theorem \ref{Theo:GGPcondi} is the possibility of applying a broader class of linear transformations. This calls for a slight adaptation of the basic definition (for regular operators) from the theory of generalized stochastic processes, which is founded on the Minlos-Bochner theorem.
\begin{definition}[Linear transformation of a generalized stochastic process]
\label{Def:LinearTransform}
Consider the following setting:
\begin{itemize}
\item a Hilbert--Gelfand space $\Spc H$ with $\Spc S \embedD \Spc H \embedD \Spc S'$;
\item a continuous linear operator $\Op T: \Spc H \to \Spc S'$ and its adjoint $\Op T^\ast: \Spc S \to \Spc H'$;
\item a generalized stochastic process $\gRand Z$ in $\Spc S'$ whose characteristic form $\ProbF_{\gRand Z}: \Spc S \to \C$
has a continuous and positive-definite extension on $\Spc H'$.
\end{itemize}
Then, the transformed process $\Op T\{\gRand Z\}$ is the generalized stochastic process in $\Spc S'$ whose
characteristic functional is $\ProbF_{\Op T\{\gRand Z\}}(\varphi)=\ProbF_{\gRand Z}(\Op T^\ast \varphi)$.

For the class of regular operators, the pair of spaces $(\Spc H,\Spc H')$ can be substituted by $(\Spc S',\Spc S)$ with $\Op T^\ast: \Spc S \to \Spc S$.
\end{definition}

In the same vein,
the application of a regular operator  $\Op L: \Spc S' \to \Spc S'$ to a generic GGP $\gRand G \sim \Spc N(\mu,\Op Q)$ yields a transformed process $\Op L\{\gRand G\} \sim \Spc N(\Lop\mu,\Op L\Op Q\Op L^\ast)$. The latter is a GGP in $\Spc S'$ with mean $\Lop\mu\in \Spc S'$ and covariance operator
$\Op L\Op Q\Op L^\ast: \Spc S \to \Spc S'$, with the transformation rules being the same as in finite dimension.

\subsection{Innovation Models and Gaussian Solutions of SDEs}
\label{Sec:GaussianSDE}
We shall now use the operator-theoretic tools of Section \ref{Sec:NativeRecap} to identify the important subclass of GGP that admit an innovation representation.
We first consider a proper GGP in $\Spc S'$ whose positive-definite covariance operator is factorizable as $\Op A=\Op T\Op T^\ast$, with $\Op T: L_2 \to  \Spc H_\Lop\embedC \Spc S'$.
Let us assume that the factor $\Op T$ is invertible, with an inverse $\Lop=\Op T^{-1}$ that is admissible with $\V p=\{0\}$ (Section \ref{Sec:NativeRecap}); the positive-definiteness of $(\Op T\Op T^\ast): \Spc S\to \Spc S'$ is then automatic, since $\langle \Op T\Op T^\ast \varphi,\varphi\rangle=\|\Op T^\ast \varphi\|^2_{L_2}$ vanishes only when $\varphi=0$.
We now show that such a process can be synthesized by applying the linear transformation  $\Op T=\Lop^{-1}$ to the Gaussian white noise (or innovation process) $\gRand w\sim \Spc N(0,\Identity)$. In our formalism, this generation mechanism is expressed as
\begin{align}
\label{Eq:GenerateG}
\gRand G(\varphi)=\Op T\{\gRand W\}(\varphi)=\langle\Op T\{\gRand W\},\varphi\rangle=
\langle\gRand W,\Op T^\ast\{\varphi\}\rangle=\gRand w(\Op T^\ast\varphi),
\end{align}
where $\Op T^\ast: (\Spc S')'=\Spc S \to (L_2)'=L_2$ is the adjoint of $\Op T: L_2 \to \Spc S'$.
The continuity of $\Op T$ together with the property that the domain of the white noise $\gRand w$ is extendable to $L_2=(L_2)'$ (see Section \ref{Sec:ExtendedProcesses}) ensure that the right-hand side of \eqref{Eq:GenerateG} is well-defined for any $\varphi \in \Spc S$.
This leads to the conclusion that
\begin{align*}
\ProbF_{\gRand G}(\varphi)=\ProbF_{\gRand w}(\Op T^\ast \varphi) 
=\exp\left(-\tfrac{1}{2} \|\Op T^\ast\varphi\|_{L_2}^2\right)=\exp\left(-\tfrac{1}{2} \langle(\Op T \Op T^\ast) \varphi, \varphi\rangle\right),
\end{align*}
which proves that $\gRand G \sim \Spc N(0, \Op A)$.
We also observe that the factorization $\Op A=\Op T\Op T^\ast$ is not unique: there is a large equivalence class of possible transformations given by
$\tilde{\Op T}=\Op T\Op U $, where $\Op U$ is an arbitrary unitary operator $L_2 \to L_2$ with the property that $\Op U^{-1}=\Op U^\ast$.

Next, we move to the scenario where the ``shaping'' operator $\Op T=\Op L^{-1}_{\V \phi}$ is the distributional extension of the operator specified by 
\eqref{III-Eq:gInverse} (right-inverse of a $\V p$-admissible
 $\Lop$). While this yields a generalized process
 $\gRand G=\Op L^{-1}_{\V \phi}\{\gRand W\}$ whose covariance operator is only positive
semi-definite, Theorem \ref{Theo:GGPcondi} guarantees that it is well-defined, with the
$\V p$-conditionally positive-definite Gram operator $\Op A_{\V \phi}$ of \eqref{III-Eq:Aphi}
as its covariance. To record that $\gRand G$ is the unique GGP that is whitened by $\Lop$ and
satisfies the boundary conditions $\V \phi(\gRand G)=\V 0$, we write its covariance operator as
\begin{align}
\Op C_{\Lop, \V \phi}=\Op L_\V \phi^{-1} \Op L_\V \phi^{-1\ast}=\Op A_\V \phi|_{\Spc S}: \Spc S \to \Spc S'.
\end{align}
and write $\gRand G \sim \Spc N(0, \Op C_{\Lop, \V \phi})$.

The $\V p$-conditional-positiveness of $\Op C_{\Lop, \V \phi}$ translates into the GGP being {\em degenerate} in the sense that it vanishes almost surely on the non-trivial null space $\Spc N_\V \phi$ of the extended operator $\Op A_\V \phi: \Spc H_\Lop' \to \Spc H_\Lop$ (see Theorem \ref{Theo:GGPcondi}).
Yet, unlike degenerate scenarios in finite dimensions, $\gRand G$ possesses an innovation $\gRand W$ that can be recovered by applying the {\em whitening operator}
$\Lop$ as in
\begin{align}
\label{Eq:Innovationmodel0}
\Lop\{\gRand G\}: \varphi \mapsto \langle\Op \Lop\{\gRand G\},\varphi\rangle=
\langle\Op \Lop\Op L^{-1}_{\V \phi}\{\gRand W\},\varphi \rangle=\gRand W(\Op L^{{-1}\ast}_{\V \phi}\, \Lop^\ast\varphi)=W(\varphi),
\end{align}
where the adjoint manipulation is legitimate because: (i) the adjoint operator $\Lop^\ast$ continuously maps $\Spc S \embedC \Spc H'_\Lop\to (L_2')=L_2$ (the extended domain of $\gRand W$); and (ii)
$\Op L^{{-1}\ast}_{\V \phi}\Lop^\ast=(\Lop\Op L^{-1}_{\V \phi})^\ast=\Identity$ on $L_2$, because $\Op L^{-1}_{\V \phi}$ is a right-inverse of $\Lop$ (Section \ref{Sec:NativeRecap}).
These properties imply that the process $\gRand g$, or rather its realizations $g \in \Spc S'$, can also be specified as 
solutions of a linear stochastic differential equation (SDE) driven by a Gaussian white noise $w \in \Spc S'$;
namely,
\begin{align}
\label{Eq:SDE}
\Lop\{\gRand G\}=\gRand W \qquad \mbox{ subject to } \qquad \V \phi(\gRand G)=\V 0,
\end{align}
which is to be understood in the weak (or distributional) sense.
The intuition there is that the ambiguity induced by the nontrivial null space of $\Lop$ gets resolved by imposing $N_0$ linear boundary conditions.

We are now in a position to state the property that ties the present paper to the
variational theory of \cite{Unser2026a,Unser2026b}: the second-order structure of $\gRand G$
and the Hilbertian structure of the native space are one and the same, up to an isometry. We recall that the {\em Gaussian
space} of a process is the closure, in quadratic mean, of the linear span of its
observations \cite{Janson1997}.
\begin{proposition}[Congruence between the native space and the Gaussian space]
\label{Theo:Congruence}
Let $\Lop$ be $\V p$-admissible and let $\gRand G \sim \Spc N(0,\Op C_{\Lop, \V \phi})$ be
the generalized Gaussian process that is whitened by $\Lop$ under the boundary conditions
$\V \phi(\gRand G)=\V 0$, with Gaussian space
\begin{align*}
\Spc G_{\gRand G}=\overline{{\rm span}}\big\{\gRand G(\varphi): \varphi \in \Spc S\big\}
\subset L_2(\Omega,\Prob).
\end{align*}
The Gaussian space is thus intrinsic to the process as initially defined on $\Spc S$
(Definition \ref{Theo:GGauss}); the same space is obtained by letting the argument range
over the whole of $\Spc H'_\Lop$, because the extended variables of
Theorem \ref{Theo:GGPcondi} are, by construction, quadratic-mean limits of variables
$\gRand G(\varphi)$ with $\varphi \in \Spc S$. Then, the following properties hold.
\begin{enumerate}
\item The map $\nu \mapsto \gRand G(\nu)$ annihilates $\Spc N_{\V \phi}$ and factors
through $\Op P_{\V \phi}$. Its restriction
\begin{align}
\label{Eq:Congruence}
\Spc U: \Spc H_{\Op A} \to \Spc G_{\gRand G},\qquad \Spc U\{\nu\}=\gRand G(\nu),
\end{align}
is a unitary isomorphism, since
$\Exp\{\gRand G(\nu)\gRand G(\upsilon)\}=\langle \Op A_{\V \phi} \nu,\upsilon\rangle=\langle \nu,\upsilon\rangle_{\Spc H'_\Lop}$
for all $\nu,\upsilon \in \Spc H_{\Op A}$.
\item By composition with the Riesz map
$\Op A_{\V \phi}: \Spc H_{\Op A} \to \Spc H_{\Lop, \V \phi}$ of
Section \ref{Sec:NativeRecap}, the primary native space $\Spc H_{\Lop, \V \phi}$
is congruent to $\Spc G_{\gRand G}$ as well, with
$\Op A_{\V \phi}\{\nu\} \mapsto \gRand G(\nu)$.
\item Both spaces are congruent to $L_2$ through the innovation, in the sense that
$\gRand G(\nu)=\gRand W(\Op L^{-1\ast}_{\V \phi}\nu)$ and
$\|\Op L^{-1\ast}_{\V \phi}\nu\|_{L_2}=\|\nu\|_{\Spc H'_\Lop}$ for every
$\nu \in \Spc H_{\Op A}$.
\end{enumerate}
\end{proposition}

\begin{proof}
\medskip
\emph{Isometry}: The covariance identity is \eqref{Eq:Cov2} of
Theorem \ref{Theo:GGPcondi}, restricted to $\Spc H_{\Op A}$ and combined with the
Riesz pairing \eqref{Eq:RieszPairingIII}. It shows that $\Spc U$ is a linear isometry from
$\Spc H_{\Op A}$ into $L_2(\Omega,\Prob)$, the action of $\gRand G$ on the whole of
$\Spc H'_\Lop$ being granted by Theorem \ref{Theo:GGPcondi}.

\medskip
\emph{Surjectivity}: Every $\varphi \in \Spc S$ decomposes as
$\varphi=\Op P_{\V \phi}\varphi + \phi$ with $\Op P_{\V \phi}\varphi \in \Spc H_{\Op A}$ and
$\phi \in \Spc N_{\V \phi}$, while $\gRand G(\phi)=0$ almost surely by
Theorem \ref{Theo:GGPcondi}. This yields
$\gRand G(\varphi)=\gRand G(\Op P_{\V \phi}\varphi)=\Spc U\{\Op P_{\V \phi}\varphi\}$, so
that the range of $\Spc U$ contains the generating family
$\{\gRand G(\varphi): \varphi \in \Spc S\}$ and is therefore dense in $\Spc G_{\gRand G}$. That range is also closed, because $\Spc U$ is an isometry,
and the two properties combine into surjectivity.

\medskip
\emph{Innovation form}: The identity
$\gRand G(\nu)=\gRand W(\Op L^{-1\ast}_{\V \phi}\nu)$ is \eqref{Eq:GenerateG} with
$\Op T=\Op L^{-1}_{\V \phi}$, while the norm identity follows from the factorization
\eqref{III-Eq:Factorization}, which gives
$\|\Op L^{-1\ast}_{\V \phi}\nu\|^2_{L_2}=\langle \Op A_{\V \phi}\nu,\nu\rangle$.
\end{proof}

Proposition \ref{Theo:Congruence} is the precise sense in which the deterministic and the
stochastic constructions coincide: the Riesz structure that \cite{Unser2026b} attaches to
$\Lop$ and the second-order structure that the present paper attaches to $\gRand G$ are
related by an isometry, while the innovation $\gRand W$ provides an orthonormal coordinate
system for both. In particular, the reproducing kernel of
\cite[Section~\ref{II-Sec:RKHS}]{Unser2026b} is the covariance function of $\gRand G$, which recovers the classical correspondence of Lo\`eve, Parzen and
Kailath \cite{Parzen1961,Kailath1970} at the level of generality of the present framework. The
pointwise form of that identity, together with the explicit kernel that the boundary
conditions produce, is Theorem \ref{Theo:ReproKernel}.

While the formal description $\gRand G=\Op L^{-1}_{\V \phi}\{\gRand W\}$ of a GGP looks rather innocuous, it is a remarkably powerful alternative to the classic integral representation of Gaussian stochastic processes. We shall make our point by considering a
famous instance of such a representation: Mandelbrot's definition \cite{Mandelbrot1968} of fractional Brownian motion (fBm) with Hurst index $H \in (0,1)$ as 
\begin{align}
\label{Eq:fBMH}
\Rand B_H(t)=\frac{1}{\Gamma(H+\frac{1}{2})}\int_{\R }\Big((t-\tau)^{H-\frac{1}{2}}_+-(-\tau)_+^{H-\frac{1}{2}}\Big)\dint \Rand B(\tau)
\end{align}
where $\Rand B(\cdot)$ is the Wiener process (considered as a random measure) and $\dint\Rand B$ the representation of Gaussian white noise in the classic notation of stochastic calculus.
The equivalent description of fBm in the present framework is
$B_H=\Op L^{-1}_{\V \phi}\{\Rand W\}$ with $\Lop= \Dop^{\gamma}$, $\gamma=H+\frac{1}{2}\in(\tfrac{1}{2},\tfrac{3}{2})$, $N_0=1$, and $(p_1,\phi_1)=(1,\delta)$. The operator $\Dop^{\gamma}$
is the fractional derivative of order $\gamma$ with Fourier symbol $(\jj \omega)^\gamma$, which induces fractional splines \cite{Unser2000}.
The relevant Green's function can be found in \cite[Table~\ref{II-Tab:Derivatives}]{Unser2026b}.
One readily verifies that the integral kernel in \eqref{Eq:fBMH} is the Schwartz kernel of the fractional integral operator $\Op L^{-1}_{\V \phi}=\Op D^{-\gamma}_{\delta}$. This equivalence also yields a simpler description of fBm as the solution of the fractional SDE
\begin{align}
\Dop^{H +\frac{1}{2}}\{\Rand B_H\}=\Rand W \quad  \mbox{ s.t.} \quad B_H(0)=0,\end{align}
which provides immediate insights on the fractal and regularity properties of these processes \cite{Blu2007a}.
As expected, and in accordance with Theorem \ref{Theo:ReproKernel}, the correlation function
of fBm coincides with the Schwartz kernel of the covariance operator: 
\begin{align}
 \label{Eq:CovfBm}
 \Exp\{\Rand B_H(t)\Rand B_H(\tau) \}=\Op A_{\V \phi}\{\delta_\tau\}(t)=\tfrac{\sigma^2_H}{2}\big(|t|^{2H}+|\tau|^{2H}-|t-\tau|^{2H} \big) \end{align}
where $\sigma^2_H=\frac{1}{\sin(\pi H)\Gamma(2H+1)}$ is the normalization that corresponds to the Green's function $\rho_\Lop(x)=\frac{x_+^{\gamma-1}}{\Gamma(\gamma)}$ of \cite[Table~\ref{II-Tab:Derivatives}]{Unser2026b},
with the case of Brownian motion being recovered for $H=\tfrac{1}{2}$ with $\Lop=\Dop$. Appendix \ref{App:fBm} carries out this computation for an arbitrary number $N_0$ of boundary conditions, which yields the family of higher-order fractional Brownian motions, while Section \ref{Sec:fBf} transposes it to $\R^d$, which yields L\'evy's Brownian field and its fractional extensions.

\section{Link with the Classical Theory}
\label{Sec:ReproKernel}
In this section, we take $(\Spc S,\Spc S')=\big(\Spc S(\R^d),\Spc S'(\R^d)\big)$ and make the link between our operator-theoretic formulation and the pointwise language of kernels, which is the one that is used in the theory of reproducing-kernel Hilbert spaces (RKHS) \cite{Aronszajn1950,Berlinet2004,Paulsen2016,Saitoh2016}, in the classical theory of Gaussian processes \cite{Lukic2001,vanderVaart2008,Rasmussen2006,Kanagawa2018}, in scattered-data approximation, where the term ``native space'' originates \cite{Schaback1999,Wendland2005,Schaback2006,Fasshauer2015}, and in kernel-based learning \cite{Scholkopf2001,Steinwart2008,Muandet2017}; the measure-theoretic side of the correspondence between kernels and integral operators is treated in \cite{Steinwart2012}. 

The basic objects of that language are kernels and the integral operators that they specify, which we define first, with the appropriate safeguards. To control growth, we introduce the space
\begin{align}
\label{Eq:SlowGrowth}
C_{\rm sg}(\R^d)\eqdef\Big\{f \in C(\R^d): \ \exists\, \alpha \ge 0 \ \mbox{ such that } \
\sup_{\V x \in \R^d}\frac{|f(\V x)|}{(1+\|\V x\|)^{\alpha}}<\infty\Big\}
\end{align}
of continuous functions of slow growth, with
$\Spc S(\R^d) \subset C_{\rm sg}(\R^d) \subset \Spc S'(\R^d)$. Likewise, a bivariate function
$a: \R^d \times \R^d \to \R$ (referred to as a kernel) is said to be of slow growth if it is
jointly continuous and if
$|a(\V x,\V y)|\le C(1+\|\V x\|)^{\alpha}(1+\|\V y\|)^{\alpha}$ for some $C>0$ and
$\alpha \ge 0$. With any such kernel, we associate the integral operator
\begin{align}
\label{Eq:KernelOp}
\Op T_a: \varphi \mapsto \Op T_a\{\varphi\}=\int_{\R^d} a(\cdot,\V y)\, \varphi(\V y)\dint \V y.
\end{align}
The stated conditions ensure that \eqref{Eq:KernelOp} well-defined on
$\Spc S(\R^d)$: the rapid decay of $\varphi$ against the polynomial bound on $a$ makes the
integral absolutely convergent for every $\V x \in \R^d$; the joint continuity of $a$,
combined with dominated convergence, yields $\Op T_a\{\varphi\} \in C_{\rm sg}(\R^d)$ (with
the same exponent $\alpha$); and the bound
$|\langle \Op T_a\{\varphi\},\psi\rangle|\le C \int_{\R^d}\int_{\R^d}
(1+\|\V x\|)^{\alpha}(1+\|\V y\|)^{\alpha}\,|\psi(\V x)|\,|\varphi(\V y)|\dint \V x \dint \V y$,
whose right-hand side is controlled by finitely many Schwartz seminorms of $\varphi$ and
$\psi$, establishes the continuity of $\Op T_a: \Spc S(\R^d) \to \Spc S'(\R^d)$. The same
bounds allow the test function in \eqref{Eq:KernelOp} to be replaced by a measure with finite
moments of order $\alpha$---in particular, by a finite atomic measure---so that the
evaluations $\Op T_a\{\delta_{\V y}\}=a(\cdot,\V y)$ make sense pointwise.

The kernels of interest to us are the two Green's functions attached to a $\V p$-admissible
operator $\Lop$. Specifically, we assume that the right inverse $\Lop^{-1}$ and the
conditional inverse $\Op A$ of $\Op B=\Lop^\ast\Lop$ (Definition \ref{Def:Admissible}) admit
the integral
representations
\begin{align}
\label{Eq:IntegralOps}
\Lop^{-1}: \varphi \mapsto \int_{\R^d} \rho_{\Lop}(\cdot,\V y) \varphi(\V y)\dint \V y,
\qquad
\Op A: \varphi \mapsto \int_{\R^d} \rho_{\Op B}(\cdot,\V y) \varphi(\V y) \dint \V y,
\end{align}
where $\rho_\Lop$ and $\rho_{\Op B}$ are the (canonical) Green's functions of $\Lop$ and of
the normal operator $\Op B=\Lop^\ast\Lop$, respectively, and where $\rho_{\Op B}$ is
symmetric with $\rho_{\Op B}(\V x,\V y)=\rho_{\Op B}(\V y,\V x)$. The representation of
$\Lop^{-1}$ is understood in the sense of the Schwartz kernel theorem, with no regularity
requirement on $\rho_\Lop$. By contrast, whenever $\rho_{\Op B}$ is of slow growth, the
conditional inverse $\Op A=\Op T_{\rho_{\Op B}}$ falls within the setting of
\eqref{Eq:KernelOp}, so that the evaluations $\Op A\{\delta_{\V y}\}=\rho_{\Op B}(\cdot,\V y)$
make sense pointwise; this is the hypothesis under which the pointwise statements of this
section hold.

The first outcome is
Theorem \ref{Theo:ReproKernel}, which identifies a single bivariate function $a_{\V \phi}$
as three objects at once: the Schwartz kernel (generalized impulse response) of the Gram
operator $\Op A_{\V \phi}$ of \eqref{III-Eq:Aphi}; the reproducing kernel
of the primary native space $\Spc H_{\Lop,\V \phi}$ of Section \ref{Sec:NativeRecap}; and the
covariance function of the generalized Gaussian process
$\gRand g \sim \Spc N(0,\Op A_{\V \phi})$ of Theorem \ref{Theo:GGPcondi}. It is therefore
the pointwise counterpart of the congruence of Proposition \ref{Theo:Congruence}, and it is
also the formula that \cite[Section~\ref{II-Sec:RKHS}]{Unser2026b} leaves open: it identifies the
reproducing kernel with the Green's function $\rho_{\Op B}$ of the normal operator in the
non-degenerate case $\Spc N_{\V p}=\{0\}$ only, whereas the boundary conditions
$\V \phi(f)=\V 0$ modify that kernel as soon as $N_0\ge 1$.

\begin{definition}[Reproducing kernel]
\label{Def:RK}
A Hilbert space $\Spc H \embedC \Spc S'(\R^d)$ is a {\em reproducing-kernel Hilbert space}
(RKHS) if it is a space of functions on which all the sampling functionals are continuous;
that is, if $\delta_{\V x_0} \in \Spc H'$ for every $\V x_0 \in \R^d$. The {\em reproducing
kernel} of such a space is the unique function $h: \R^d \times \R^d \to \R$ such that
\begin{align}
&\mbox{(i) } h(\cdot,\V x_0) \in \Spc H \mbox{ for all }\V x_0 \in \R^d;
\nonumber\\
&\mbox{(ii) } f(\V x_0)=\langle h(\cdot,\V x_0), f\rangle_{\Spc H} \mbox{ for all }
f \in \Spc H \mbox{ and } \V x_0 \in \R^d.
\label{Eq:reproform}
\end{align}
\end{definition}
The two properties in \eqref{Eq:reproform} say that $h(\cdot,\V x_0)$ is the Riesz
representative of $\delta_{\V x_0}$, which is why the existence of the kernel and the
continuity of the sampling functionals are equivalent.

\subsection{Reproducing Kernel and Covariance Function}
We now provide the pointwise transcription of the two descriptions of the native space that
\cite{Unser2026b} provides: the conditioned parametrization
$\Op A_{\V \phi}=\Op P^\ast_{\V \phi}\Op A\Op P_{\V \phi}$ and the factorization
$\Op A_{\V \phi}=\Op L^{-1}_{\V \phi}\Op L^{-1\ast}_{\V \phi}$
\big(Eqs.~\eqref{III-Eq:Aphi} and \eqref{III-Eq:Factorization}\big). The first one
produces the explicit formula, while the second one
exhibits the kernel as a Gram (or autocorrelation) integral.

\begin{theorem}[Reproducing kernel of the primary native space]
\label{Theo:ReproKernel}
Let $\Lop: \Spc S(\R^d) \to L_2(\R^d)$ be $\V p$-admissible in the sense of
Section \ref{Sec:NativeRecap}, and let $(\V p, \V \phi)$ be a biorthogonal system with
$\V p \subset C_{\rm sg}(\R^d)$ and $\V \phi \subset \Spc H'_\Lop$, as in Section \ref{Sec:NativeRecap}, such that
the boundary functionals can be paired with the Green's functions; that is, such that
$\langle \phi_n,\rho_\Lop(\cdot,\V y)\rangle$, $\langle \phi_n,\rho_{\Op B}(\cdot,\V y)\rangle$,
and $\langle \phi_m,\Op A\{\phi_n\}\rangle$ are well-defined for all $\V y \in \R^d$.
Then, the following properties hold, Items \ref{Item:Gram}--\ref{Item:Covariance} under the
additional hypothesis that the Green's function $\rho_{\Op B}$ of \eqref{Eq:IntegralOps} is of
slow growth.
\begin{enumerate}
\item \label{Item:KernelInverse} The conditioned inverse
$\Op L^{-1}_{\V \phi}=\Op P^\ast_{\V \phi}\Lop^{-1}: L_2(\R^d) \to \Spc H_{\Lop,\V \phi}$ \big(see \eqref{III-Eq:gInverse}\big) is the integral operator with kernel
\begin{align}
\label{Eq:kernelinverseg}
g_{\V \phi}(\V x,\V y) =\rho_{\Lop}(\V x,\V y)- \sum_{n=1}^{N_0}  p_n(\V x) q_n(\V y),
\qquad
q_n=\Lop^{-1\ast}\{\phi_n\},
\end{align}
where the functions $q_n(\V y)=\langle \phi_n,\rho_{\Lop}(\cdot,\V y)\rangle$ are the images of
the boundary functionals under the adjoint of the Green's operator. No regularity of
$\rho_\Lop$ is required here: \eqref{Eq:kernelinverseg} is an identity between Schwartz
kernels, which is what allows it to cover the leading examples $\Lop=\Dop$ and
$\Lop=\Dop^\gamma$ with $\gamma<1$, whose Green's functions are discontinuous at the diagonal.
\item \label{Item:Gram} The Gram operator $\Op A_{\V \phi}$ is the integral operator with the
symmetric kernel
\begin{align}
\label{Eq:kernelrepro}
a_{\V \phi}(\V x,\V y)
&=\rho_{\Op B}(\V x,\V y)- \sum_{n=1}^{N_0}  p_n(\V x) v_n(\V y)-
\sum_{n=1}^{N_0}  v_n(\V x) p_n(\V y) \nonumber \\
&\quad + \sum_{m=1}^{N_0}\sum_{n=1}^{N_0}r_{m,n}\,p_m(\V x)p_n(\V y),
\end{align}
where
\begin{align}
\label{Eq:vn}
v_n&=\Op A\{\phi_n\}, \qquad \mbox{i.e., } \
v_n(\V y)=\int_{\R^d} \phi_n(\V z)\rho_{\Op B}(\V z, \V y)  \dint \V z,\\
\label{Eq:rmn}
r_{m,n}
&=\langle \phi_m, \Op A\{\phi_n\} \rangle=
\int_{\R^d}\int_{\R^d} \phi_m(\V x)\rho_{\Op B}(\V x,\V y)\phi_n(\V y) \dint \V x\dint \V y.
\end{align}
\item \label{Item:RKHS} The three following properties are equivalent: (a)
$\delta_{\V x} \in \Spc H'_\Lop$ for every $\V x \in \R^d$; (b)
$g_{\V \phi}(\V x,\cdot) \in L_2(\R^d)$ for every $\V x \in \R^d$; and (c) the primary native
space $\Spc H_{\Lop,\V \phi}$ is an RKHS. When they hold, the reproducing kernel of
$\Spc H_{\Lop,\V \phi}$ is $a_{\V \phi}$, which is then also the autocorrelation of the kernel
\eqref{Eq:kernelinverseg},
\begin{align}
\label{Eq:kernelfactor}
a_{\V \phi}(\V x,\V y)=\int_{\R^d} g_{\V \phi}(\V x, \V z) g_{\V \phi}(\V y, \V z)\dint \V z
=\langle g_{\V \phi}(\V x,\cdot),g_{\V \phi}(\V y,\cdot)\rangle_{L_2},
\end{align}
with $a_{\V \phi}(\V x,\V x)=\|g_{\V \phi}(\V x,\cdot)\|^2_{L_2}$; and the reproducing kernel of
the full native space $\Spc H_\Lop=\Spc H_{\Lop,\V \phi}\oplus \Spc N_{\V p}$ for the inner
product \eqref{III-Eq:InnerProd} is
$a_{\V \phi}(\V x,\V y)+\sum_{n=1}^{N_0}p_n(\V x)p_n(\V y)$.
\item \label{Item:Covariance} Under the conditions of Item \ref{Item:RKHS}, the process
$\gRand g \sim \Spc N(0,\Op A_{\V \phi})$ of Theorem \ref{Theo:GGPcondi} admits the
pointwise reading $\gRand g(\V x)\eqdef\gRand g(\delta_{\V x})$, and $a_{\V \phi}$ is its
covariance function:
\begin{align}
\label{Eq:CovarianceFunction}
\Exp\big\{\gRand g(\V x)\gRand g(\V y)\big\}=a_{\V \phi}(\V x,\V y),
\qquad (\V x,\V y) \in \R^d \times \R^d.
\end{align}
\end{enumerate}
\end{theorem}
The pairing condition on $\V \phi$ is automatic for $\V \phi \subset \Spc S(\R^d)$, and it
holds for the Dirac-type systems $\phi_n=\delta^{(n-1)}$ and $\phi_{\V k}=\partial^{\V k}\delta$
of Appendix \ref{App:fBm} and of Section \ref{Sec:fBf}, whose action amounts to the point evaluation
of the (sufficiently smooth) kernels and of their derivatives.

\begin{proof}
\medskip
\emph{Kernel of the conditioned inverse (Item \ref{Item:KernelInverse})}:
By \cite[Eq.~\eqref{II-Eq:gInverse}]{Unser2026b}, we have that
$\Op L^{-1}_{\V \phi}=\Op P^\ast_{\V \phi}\Op G^\ast=(\Identity-\Op R_{\V p}\Op R_{\V \phi})\Lop^{-1}$,
so that $\Op P^\ast_{\V \phi}\{f\}=f-\sum_{n=1}^{N_0} p_n\langle \phi_n,f\rangle$ acts on the
first argument of the kernel $\rho_\Lop$ of $\Lop^{-1}$. Since
$\langle \phi_n,\rho_\Lop(\cdot,\V y)\rangle=q_n(\V y)$ by the definition of $q_n$, this yields
\eqref{Eq:kernelinverseg}. The identity $q_n=\Lop^{-1\ast}\{\phi_n\}$ follows from the fact
that the Schwartz kernel of an adjoint operator is the transpose of the original one.

\medskip
\emph{Kernel of the Gram operator (Item \ref{Item:Gram})}:
We evaluate $a_{\V \phi}(\cdot,\V y)=\Op A_{\V \phi}\{\delta_{\V y}\}$ by applying the three
factors of $\Op A_{\V \phi}=\Op P^\ast_{\V \phi}\Op A\Op P_{\V \phi}$ in succession. The
projector $\Op P_{\V \phi}=\Identity-\Op R_{\V \phi}\Op R_{\V p}$ gives
$\Op P_{\V \phi}\{\delta_{\V y}\}=\delta_{\V y}-\sum_{n=1}^{N_0} p_n(\V y)\phi_n$, whose image
under $\Op A$ is $\rho_{\Op B}(\cdot,\V y)-\sum_{n=1}^{N_0} p_n(\V y) v_n$ with $v_n$ given by
\eqref{Eq:vn}. The symmetry of $\rho_{\Op B}$ yields
$\langle \phi_m,\rho_{\Op B}(\cdot,\V y)\rangle=v_m(\V y)$, while the definition
\eqref{Eq:rmn} of $r_{m,n}$ gives $\langle \phi_m,v_n\rangle=r_{m,n}$. The final application
of $\Op P^\ast_{\V \phi}$ therefore produces
\begin{align*}
a_{\V \phi}(\cdot,\V y)&=\Big(\rho_{\Op B}(\cdot,\V y)-\sum_{n=1}^{N_0} p_n(\V y) v_n\Big)
-\sum_{m=1}^{N_0} p_m \Big(v_m(\V y)-\sum_{n=1}^{N_0}r_{m,n}p_n(\V y)\Big),
\end{align*}
which is \eqref{Eq:kernelrepro}. That expression is symmetric in $(\V x,\V y)$ because
$\rho_{\Op B}$ is symmetric and because $r_{m,n}=r_{n,m}$, the operator $\Op A$ being
symmetric.

\medskip
\emph{Reproduction property (Item \ref{Item:RKHS})}:
Let us fix $\V x \in \R^d$. Since $\Spc N_{\V \phi}\subset \Spc H'_\Lop$,
the decomposition $\delta_{\V x}=\Op P_{\V \phi}\{\delta_{\V x}\}+\sum_{n}p_n(\V x)\phi_n$
shows that (a) is equivalent to
$\Op P_{\V \phi}\{\delta_{\V x}\} \in \Spc H_{\Op A}=(\Spc H_{\Lop,\V \phi})'$, which is
precisely the continuity of the sampling functional on $\Spc H_{\Lop, \V \phi}$; that is,
(a) $\Leftrightarrow$ (c) by Definition \ref{Def:RK}. To bring in (b), we parametrize
$\Spc H_{\Lop,\V \phi}=\Op L^{-1}_{\V \phi}(L_2)$ by the innovation. Every
$f \in \Spc H_{\Lop,\V \phi}$ is $f=\Op L^{-1}_{\V \phi}\{w\}$ with
$w=\Lop f \in L_2(\R^d)$ and $\|f\|_{\Spc H_{\Lop,\V \phi}}=\|w\|_{L_2}$
\big(by \eqref{III-Eq:InnerProd}, since $\V \phi(f)=\V 0$\big), so that
Item \ref{Item:KernelInverse} gives
$f(\V x)=\int_{\R^d} g_{\V \phi}(\V x,\V z)w(\V z)\dint \V z$. The functional
$w \mapsto f(\V x)$ is therefore bounded on $L_2(\R^d)$ if and only if
$g_{\V \phi}(\V x,\cdot) \in L_2(\R^d)$, by the Riesz representation theorem, which is the
equivalence of (b) and (c). Let us now assume that (a) holds and let
$f \in \Spc H_{\Lop,\V \phi}$, so that $\V \phi(f)=\V 0$. The
biorthogonality then gives
$\langle \Op P_{\V \phi}\{\delta_{\V x}\},f\rangle=f(\V x)-\sum_{n=1}^{N_0}p_n(\V x)\langle \phi_n,f\rangle=f(\V x)$,
where $\Op P_{\V \phi}\{\delta_{\V x}\} \in \Spc H_{\Op A}$. Since
$\Op A_{\V \phi}: \Spc H_{\Op A} \to \Spc H_{\Lop,\V \phi}$ is the Riesz map of Section \ref{Sec:NativeRecap}
\big(see \eqref{Eq:RieszPairingIII}\big), and since
$\Op A_{\V \phi}\Op P_{\V \phi}=\Op A_{\V \phi}$ by idempotence of the projector, this
translates into
\begin{align*}
f(\V x)=\langle \Op A_{\V \phi}\Op P_{\V \phi}\{\delta_{\V x}\},f\rangle_{\Spc H_{\Lop,\V \phi}}
=\langle a_{\V \phi}(\cdot,\V x),f\rangle_{\Spc H_{\Lop,\V \phi}},
\end{align*}
which is the reproduction property \eqref{Eq:reproform}, the membership
$a_{\V \phi}(\cdot,\V x) \in \Spc H_{\Lop,\V \phi}$ being granted by the range of
$\Op A_{\V \phi}$. This proves that (a) implies (c), while the converse implication is
Definition \ref{Def:RK}. The formula for the reproducing kernel of $\Spc H_\Lop$ is obtained in the same way from
the Riesz map $\Op J_{\Spc H'_\Lop}=\Op A_{\V \phi}+\Op R_{\V p}$ of the full space
\big(the augmented operator $\Op A_+$ of Section \ref{Sec:NativeRecap}\big), whose second
block contributes
$\Op R_{\V p}\{\delta_{\V y}\}=\sum_{n=1}^{N_0} p_n\, p_n(\V y)$. Finally, the identity
\eqref{Eq:kernelfactor} is the pointwise transcription of the factorization \eqref{III-Eq:Factorization}: the kernel of a composition is the contraction of the two
kernels, that of $\Op L^{-1\ast}_{\V \phi}$ being the transpose of
\eqref{Eq:kernelinverseg}, and the contraction converges precisely because
$g_{\V \phi}(\V x,\cdot) \in L_2(\R^d)$.

\medskip
\emph{Covariance function (Item \ref{Item:Covariance})}:
By Theorem \ref{Theo:GGPcondi}, the action of $\gRand g$ extends to any test function in
$\Spc H'_\Lop$, which under (a) includes the sampling functionals $\delta_{\V x}$; the
pointwise reading is therefore legitimate. Equation \eqref{Eq:Cov2} with
$\varphi_1=\delta_{\V x}$ and $\varphi_2=\delta_{\V y}$ then gives
$\Exp\{\gRand g(\V x)\gRand g(\V y)\}=\langle \Op A_{\V \phi}\{\delta_{\V x}\},\delta_{\V y}\rangle=a_{\V \phi}(\V y,\V x)$,
which is \eqref{Eq:CovarianceFunction} by symmetry.
\end{proof}

Three comments are in order.

First, Formula \eqref{Eq:kernelrepro} is a double Taylor-type remainder: the raw Green's
function $\rho_{\Op B}$ is corrected in each variable by the finite-rank term that enforces the
boundary conditions, and the double sum restores the symmetry. The characteristic property of
$a_{\V \phi}$ is that $\V \phi\big(a_{\V \phi}(\cdot,\V y)\big)=\V 0$ for every
$\V y \in \R^d$, which is the pointwise transcription of the boundary conditions
$\V \phi(\gRand G)=\V 0$ of the stochastic differential equation \eqref{Eq:SDE}. For
$\Spc N_{\V p}=\{0\}$, all correction terms disappear and one recovers
$a_{\V \phi}=\rho_{\Op B}$, which is the identification of \cite[Section~\ref{II-Sec:RKHS}]{Unser2026b}
with $\rho_{\Op B}=h$; in the shift-invariant case, this is the kernel of
\cite[Eq.~\eqref{II-Eq:LSIKernel}]{Unser2026b}.

Second, the kernel $g_{\V \phi}$ of \eqref{Eq:kernelinverseg} is the impulse response of the
innovation model of Section \ref{Sec:GaussianSDE}: once $\gRand G=\Op L^{-1}_{\V \phi}\{\gRand W\}$
is available, \eqref{Eq:kernelinverseg} is the integral representation
$\gRand G(\V x)=\int_{\R^d} g_{\V \phi}(\V x,\V z)\gRand W(\V z)\dint \V z$ of the process, and
\eqref{Eq:kernelfactor} is the classical formula that expresses a covariance as the
autocorrelation of a shaping kernel. The worked example is fractional Brownian motion:
$g_{\V \phi}$ is the kernel \eqref{Eq:fBmKernel} of Appendix \ref{App:fBm}, which reduces to
the integral kernel of Mandelbrot's representation \eqref{Eq:fBMH} for $N_0=1$, while
$a_{\V \phi}$ is the covariance function \eqref{Eq:fBmCovGeneral}, which reduces to
\eqref{Eq:CovfBm}.

Third, the requirement $\delta_{\V x} \in \Spc H'_\Lop$ is not a technicality but the exact
frontier between the two readings of the theory. It is simultaneously the condition under
which the native space is a space of functions rather than of distributions, and the condition
under which the generalized process has samples $\gRand g(\V x)$ rather than only smoothed
observations $\gRand g(\varphi)$. Gaussian white noise, for which $\Spc H'_\Lop=L_2(\R^d)$ and
$\delta_{\V x} \notin L_2(\R^d)$, is the canonical counterexample on both counts.

\subsection{Conditional Positive Definiteness}
\label{Sec:Conditional}
The price to pay for a regularization operator $\Lop$ whose null space is non-trivial is that
the kernel $\rho_{\Op B}$ of the conditional inverse $\Op A=(\Lop^\ast\Lop)^{-1}$ cannot be
positive-definite. Instead, one has to settle for a weaker form of the property that factors
out the components that compromise uniqueness. The operator-level version of that notion is
the $\V p$-conditional positive definiteness of Section \ref{Sec:NativeRecap}. We now give its kernel-level counterpart, which is the form in which the property appears in the literature on radial basis
functions and conditionally positive-definite kernels
\cite{Micchelli1986,Buhmann2003,Wendland2005,Iske2018}.

\begin{definition}[Conditionally positive-definite kernel]
\label{Def:CondPositiveDef1}
Let $\V p=(p_1,\dots,p_{N_0})$ be a linearly independent family in $C_{\rm sg}(\R^d)$, with
$\Spc N_{\V p}={\rm span}\{p_n\}_{n=1}^{N_0}$ the subspace that it spans. Then, a symmetric kernel
$h: \R^d \times \R^d \to \R$ is said to be {\em $\V p$-conditionally positive-semi-definite}
if
\begin{align}
\label{Eq:CondPDkernel}
\sum_{m=1}^N \sum_{n=1}^N z_m h(\V x_m, \V x_n) z_n \ge  0
\end{align}
for any $N \in \N$, any $\V x_1, \dots, \V x_N \in \R^d$ and any $z_1, \dots, z_N \in \R$ that
are subject to the moment conditions
\begin{align}
\label{Eq:MomentCond}
\sum_{m=1}^N z_m p_n(\V x_m)=  0, \qquad n=1,\dots,N_0.
\end{align}
It is said to be {\em $\V p$-conditionally positive-definite} if the inequality
\eqref{Eq:CondPDkernel} is strict for all distinct $\V x_1, \dots, \V x_N \in \R^d$ and all
$(z_1, \dots, z_N) \in \R^N\backslash\{\V 0\}$ that satisfy \eqref{Eq:MomentCond}. The
unconditional versions of the two properties are obtained for $N_0=0$, in which case
\eqref{Eq:MomentCond} is void.
\end{definition}
The probability literature drops the qualifier ``semi'' in \eqref{Eq:CondPDkernel}. For a
shift-invariant kernel
$h(\V x,\V y)=\rho(\V x-\V y)$ with $N_0=0$, \eqref{Eq:CondPDkernel} is exactly the Bochner
condition \eqref{Eq:PSfunction} of Theorem \ref{Theo:Minlos}, applied to $\rho$ with $\R^d$
in place of $\Spc S$; this special case is the one that underlies the classical RKHS theory.

The two levels---operators and kernels---are equivalent as soon as the kernel is a genuine
bivariate function of slow growth.

\begin{theorem}[Kernel of a positive-semi-definite operator]
\label{Theo:OptoKernelPosDef}
Let $\Op A=\Op T_a$ be a symmetric integral operator of the form \eqref{Eq:KernelOp} whose
kernel $a$ is of slow growth, and let $\V p \subset C_{\rm sg}(\R^d)$. Then,
the ($\V p$-conditional) positive semi-definiteness of the operator
$\Op A: \Spc S(\R^d) \to \Spc S'(\R^d)$ in the sense of
Section \ref{Sec:NativeRecap} is equivalent to the ($\V p$-conditional)
positive semi-definiteness of its kernel $a: \R^d \times \R^d \to \R$ in the sense of
Definition \ref{Def:CondPositiveDef1}.
\end{theorem}

\begin{proof}
The two implications rest on two complementary discretizations: mollified Dirac combs for the
direct part, and Riemann sums for the converse.

\medskip
\emph{Mollified Dirac combs}:
For $k \in \N^+$, let
$u_k(\V x)=\left(\frac{k}{\sqrt{2\pi}}\right)^{d}\ee^{-\frac{1}{2}\|k \V x\|^2}$ be the
rescaled and renormalized Gaussian pulse, which is in $\Spc S(\R^d)$ with
$\langle u_k, 1\rangle=1$ and which converges to the Dirac distribution in the weak-$\ast$
topology of $\Spc M(\R^d)$ as $k \to \infty$. Given $\V x_1, \dots,\V x_N \in \R^d$ and
$z_1,\dots,z_N \in \R$, we set $\varphi_k=\sum_{n=1}^N z_n u_k(\cdot-\V x_n) \in \Spc S(\R^d)$.
The family $u_k(\cdot-\V x_m)\otimes u_k(\cdot-\V x_n)$ is an approximate identity at
$(\V x_m,\V x_n)$, so that, for a given $\epsilon>0$, we may split the corresponding integral
of $a$ into the contribution of a ball around $(\V x_m,\V x_n)$ on which the continuity of $a$
confines the integrand to $a(\V x_m,\V x_n)\pm \epsilon$, and the contribution of the
complement, over which the Gaussian decay of $u_k$ dominates the polynomial bound on $a$ and
which therefore vanishes as $k \to \infty$. This yields
\begin{align}
\label{Eq:samplingrel}
\lim_{k\to \infty} \langle \Op A \{\varphi_k\},\varphi_k\rangle
&=\lim_{k\to \infty} \int_{\R^d \times \R^d} a(\V x,\V y)(\varphi_k\otimes \varphi_k)(\V x,\V y)
\;\dint (\V x,\V y)\nonumber \\
&=\sum_{m=1}^{N} \sum_{n=1}^{N} z_m z_n a(\V x_m,\V x_n).
\end{align}

\medskip
\emph{Operator $\Rightarrow$ kernel}:
In the unconditional case, the left-hand side of \eqref{Eq:samplingrel} is non-negative for
every $k$, so that the limit is non-negative as well, which is \eqref{Eq:CondPDkernel}.
In the conditional case, we replace $\varphi_k$ by
$\Op P_{\V \phi}\{\varphi_k\}=\varphi_k-\sum_{n=1}^{N_0} c_{n,k}\phi_n \in \Spc S_{\V p}(\R^d)$
with $c_{n,k}=\langle p_n,\varphi_k\rangle$, where $\V \phi \subset \Spc S(\R^d)$ is
biorthogonal to $\V p$ (Section \ref{Sec:NativeRecap}). The positivity hypothesis
now applies, since $\Op P_{\V \phi}\{\varphi_k\} \in \Spc S_{\V p}(\R^d)$ by construction. The
coefficients satisfy $c_{n,k} \to \sum_{m=1}^N z_m p_n(\V x_m)=0$, the moment conditions
\eqref{Eq:MomentCond} being in force, so that the correction terms
$-2\sum_n c_{n,k}\langle \Op A \phi_n,\varphi_k\rangle+\sum_m\sum_n c_{m,k}c_{n,k}\langle \Op A\phi_m,\phi_n\rangle$
vanish as $k \to \infty$. The limit of
$\langle \Op A \Op P_{\V \phi}\varphi_k,\Op P_{\V \phi}\varphi_k\rangle \ge 0$ is therefore
again the right-hand side of \eqref{Eq:samplingrel}.

\medskip
\emph{Kernel $\Rightarrow$ operator}:
Let $\varphi \in \Spc S(\R^d)$, $R>0$ and $k \in \N^+$, and let us introduce the Riemann sums
\begin{align}
\label{Eq:discretizePD}
J_{R,k}(\varphi)&=\sum_{\V m\in [-k,k]^d} \sum_{\V n\in [-k,k]^d}
z_{\V m}\, a(\V x_{\V m},\V x_{\V n})\, z_{\V n},
\qquad z_{\V m}=\left(\tfrac{R }{k}\right)^{d}\varphi(\V x_{\V m}),
\end{align}
with $\V x_{\V m}=\tfrac{\V m R}{k} \in \R^d$ and $\V m \in \Z^d$. Let us set $R=i$ and
$k=i^2$, so that the mesh $\tfrac{1}{i}$ vanishes while the domain fills $\R^d$. The
integrand $F(\V x,\V y)=\varphi(\V x)a(\V x,\V y)\varphi(\V y)$ is continuous, and the rapid
decay of $\varphi$ against the polynomial bound on $a$ gives
$|F(\V x,\V y)|\le C_N(1+\|\V x\|)^{-N}(1+\|\V y\|)^{-N}$ for every $N$. Given $\epsilon>0$,
we may therefore pick $T>0$ such that both the integral of $|F|$ and the Riemann sums of $|F|$
over the complement of $[-T,T]^{2d}$ are below $\epsilon$, uniformly in $i$ (the mesh being at
most $1$); on the compact part, the uniform continuity of $F$ makes the sums converge to the
integral. Hence
$\lim_{i\to \infty}J_{i,i^2}(\varphi)=\langle \Op A\varphi,\varphi\rangle$. In the
unconditional case, $J_{i,i^2}(\varphi)\ge 0$ by \eqref{Eq:CondPDkernel}, and the conclusion
follows.

In the conditional case, the weights $z_{\V m}$ of \eqref{Eq:discretizePD} do not fulfil the
moment conditions \eqref{Eq:MomentCond} exactly, but only asymptotically, so that a correction
is needed. Let us fix $N_0$ points $\V y_1,\dots,\V y_{N_0} \in \R^d$ such that the matrix
$\M P=[p_n(\V y_j)] \in \R^{N_0\times N_0}$ is invertible, which is possible because the $p_n$
are linearly independent continuous functions. For $\varphi \in \Spc S_{\V p}(\R^d)$, we then
form the residuals $\epsilon_{i,n}=\sum_{\V m} z_{\V m} p_n(\V x_{\V m})$ and the corrective
weights $\M c_i=-\M P^{-1}\V \epsilon_i$, and we append the atoms $\V y_j$ with weights
$c_{i,j}$ to the grid. The augmented family fulfils \eqref{Eq:MomentCond} exactly, so that
Definition \ref{Def:CondPositiveDef1} applies and yields
\begin{align*}
J_{i,i^2}(\varphi)+2\sum_{j=1}^{N_0} c_{i,j}\sum_{\V m} z_{\V m}a(\V y_j,\V x_{\V m})
+\sum_{j=1}^{N_0}\sum_{j'=1}^{N_0} c_{i,j}c_{i,j'}a(\V y_j,\V y_{j'})\ge 0.
\end{align*}
Now, $\epsilon_{i,n} \to \langle p_n,\varphi\rangle=0$ as $i \to \infty$, which gives
$\M c_i \to \V 0$; moreover, the sums $\sum_{\V m} |z_{\V m}| \,|a(\V y_j,\V x_{\V m})|$ are
bounded uniformly in $i$, because $a$ is of slow growth and $\varphi$ is rapidly decaying.
The two corrective terms therefore vanish in the limit, which leaves
$\langle \Op A\varphi,\varphi\rangle \ge 0$ for all $\varphi \in \Spc S_{\V p}(\R^d)$.
\end{proof}

Neither of the two implications survives the passage to the strict versions of the
properties, and it is worth saying why, since the strict form is the one that governs the
invertibility of the interpolation matrices.

In the direction ``operator $\Rightarrow$ kernel'', the obstruction is that a strict
inequality is not preserved by the limit \eqref{Eq:samplingrel}. A counterexample on $\R$ is
the kernel $a(x,y)=x\,y\,\ee^{-(x-y)^2}$, which is of slow growth and which satisfies
$\langle \Op A \varphi,\varphi\rangle=\langle \Op A_{\rm g} \psi,\psi\rangle>0$ for every
$\varphi \in \Spc S(\R)\backslash\{0\}$, where $\psi(x)=x\varphi(x) \in \Spc S(\R)$ and where
$\Op A_{\rm g}$ is the (strictly positive-definite) Gaussian operator with kernel
$\ee^{-(x-y)^2}$. The kernel of $\Op A$, by contrast, degenerates at the origin, since $a(0,0)=0$.

In the direction ``kernel $\Rightarrow$ operator'', the obstruction is that the strict
positivity of a kernel only constrains the finite Gram matrices. The illustration is the
one-dimensional shift-invariant kernel $h_{\rm sinc}(x,y)={\rm sinc}(x-y)$, which is
positive-definite in the sense of Definition \ref{Def:CondPositiveDef1} because the functions
$\omega \mapsto \ee^{\jj \omega x_m}$ are linearly independent in $L_2([-\pi,\pi])$ for
distinct $x_m$. Its operator, by contrast, is only positive-semi-definite. Indeed, we have
that $\Fourier\{{\rm sinc}\}(\omega)={\rm rect}\big(\omega/(2 \pi)\big)$, so that
\begin{align}
\label{Eq:sincQuadratic}
\langle \Op A_{\rm sinc} \varphi,\varphi\rangle
=\frac{1}{2\pi}\int_{\R} {\rm rect}\big(\tfrac{\omega}{2\pi}\big) |\hat \varphi(\omega)|^2
\dint \omega,
\end{align}
which vanishes for any $\varphi \in \Spc S(\R)\backslash\{0\}$ whose spectrum is supported
outside the band $[-\pi,\pi]$.

What does hold, and what the theory of Section \ref{Sec:ReproKernel} actually needs, is an
exact characterization of the strict property for the kernel $a_{\V \phi}$ of
Theorem \ref{Theo:ReproKernel} in terms of the sampling functionals. Let us assume that the
conditions of Item \ref{Item:RKHS} are met and let $\V x_1,\dots,\V x_M \in \R^d$ be distinct.
For any $\V z \in \R^M$ that fulfils the moment conditions \eqref{Eq:MomentCond}, we have that
$\Op P_{\V \phi}\big\{\sum_m z_m \delta_{\V x_m}\big\}=\sum_m z_m \delta_{\V x_m}$, so that
\begin{align}
\label{Eq:GramNorm}
\sum_{m=1}^M \sum_{n=1}^M z_m\, a_{\V \phi}(\V x_m,\V x_n)\, z_n
=\Big\|\sum_{m=1}^M z_m \delta_{\V x_m}\Big\|^2_{\Spc H'_\Lop}
\end{align}
by \eqref{Eq:RieszPairingIII}. The right-hand side vanishes for some $\V z \neq \V 0$ if and
only if the $\delta_{\V x_m}$ are linearly dependent in $\Spc H'_\Lop$---the moment conditions
being automatic for such a $\V z$, as one sees by testing against $\V p \subset \Spc H_\Lop$.
We have therefore established the following.
\begin{corollary}[Strict positive definiteness at a set of nodes]
\label{Theo:StrictNodes}
Under the conditions of Item \ref{Item:RKHS} of Theorem \ref{Theo:ReproKernel}, the kernel
$a_{\V \phi}$ is strictly $\V p$-conditionally positive-definite in the sense of
Definition \ref{Def:CondPositiveDef1} {\em if and only if} the sampling functionals
$\delta_{\V x_1},\dots,\delta_{\V x_M}$ are linearly independent in $\Spc H'_\Lop$ for every
finite set of distinct points $\V x_1,\dots,\V x_M \in \R^d$.
\end{corollary}
This is the linear independence of the sampling functionals in $\Spc H'_\Lop$, which is
Condition 1 of the $\Lop$-compatibility used in Section \ref{Sec:MMSESplines} (read for
$\nu_m=\delta_{\V x_m}$). It is also the pointwise instance of the criterion of
Theorem \ref{Theo:GGPcondi}, which states that the covariance matrix of
$\big(\gRand g(\varphi_1),\dots,\gRand g(\varphi_N)\big)$ is positive-definite exactly when the
$\Op P_{\V \phi}\{\varphi_n\}$ are linearly independent in $\Spc H_{\Op A}$. The
shift-invariant case is the reassuring one: if $\widehat A \in L_1(\R^d)$ with
$\widehat A>0$ almost everywhere, then
$\sum_m\sum_n z_m z_n \rho(\V x_m-\V x_n)=\frac{1}{(2\pi)^d}\int_{\R^d}\widehat A(\bw)
\big|\sum_m z_m \ee^{-\jj \langle \bw,\V x_m\rangle}\big|^2 \dint \bw$ is strictly positive,
because the exponential sum is real-analytic and therefore vanishes only on a set of measure
zero.

\subsection{Extended Fractional Brownian Fields on $\R^d$}
\label{Sec:fBf}
We now illustrate the framework with the fractional Laplacian \cite{Kwasnicki2017,DiNezza2012}, which is the multidimensional counterpart of the fractional derivative of Section \ref{Sec:GaussianSDE}. The resulting
processes are the isotropic {\em fractional Brownian fields}: L\'evy's multiparameter
Brownian motion \cite{Levy1948} and its fractional extensions \cite{Tafti2010}. Since all
the tools are already in place, every formula below is an instance of
Theorem \ref{Theo:ReproKernel}, whose hypothesis $\V \phi \subset \Spc H'_\Lop$ covers the
Dirac-type boundary functionals.

The fractional Laplacian
$\Lop=(-\Delta)^{\gamma/2}: \Spc S(\R^d) \to \Spc S'(\R^d)$ is the shift-invariant operator whose Fourier symbol is
$\widehat{L}(\bw)=\|\bw\|^{\gamma}$, over the range
\begin{align}
\label{Eq:fBfRange}
\gamma>\tfrac{d}{2}, \qquad H=\gamma-\tfrac{d}{2} \notin \N,
\qquad k_\gamma=\lfloor H\rfloor,
\end{align}
where $H$ is the Hurst-type exponent of the resulting field and where the excluded integer
values of $H$ are the logarithmic (thin-plate) cases discussed at the end of this section.
The operator is self-adjoint, with normal operator $\Op B=(-\Delta)^{\gamma}$ of symbol
$\|\bw\|^{2\gamma}$. The fractional Laplacian meets the conditions of the multiplier criterion
\cite[Proposition~\ref{I-Prop:MultiplierConstruction}]{Unser2026a}, so that $\Lop$ is
$\V p$-admissible in the sense of Definition \ref{Def:Admissible}, with attainable null
space $\Spc N_{\Op B}=\Pi_{k_\gamma}$, the space of polynomials of degree at most
$k_\gamma$, whose dimension is $N_0=\binom{k_\gamma+d}{d}$. With the usual multi-index notation
($|\V k|=k_1+\cdots+k_d$, $\V k!=k_1!\cdots k_d!$, and
$\V x^{\V k}=x_1^{k_1}\cdots x_d^{k_d}$), the biorthogonal system is
\begin{align}
\label{Eq:fBfBiortho}
p_{\V k}(\V x)=\frac{\V x^{\V k}}{\V k!},
\qquad
\phi_{\V k}=\partial^{\V k}\delta,
\qquad |\V k|\le k_\gamma,
\end{align}
with the convention $\langle \partial^{\V k}\delta,f\rangle=\partial^{\V k}f(\V 0)$, so that the boundary conditions read $\partial^{\V k}g(\V 0)=0$ for
$|\V k|\le k_\gamma$. The primary example is the range $k_\gamma=0$, that is,
$H \in (0,1)$: the null space then reduces to the constants, with
$(p_1,\phi_1)=(1,\delta)$ and the single boundary condition $g(\V 0)=0$---L\'evy's field,
and the multidimensional analogue of Mandelbrot's normalization.

Both Green's functions are Riesz kernels. With the Fourier convention
$\hat \varphi(\bw)=\int_{\R^d}\varphi(\V x)\,\ee^{-\jj \langle \bw,\V x\rangle}\dint \V x$,
the classical transform pair of homogeneous distributions reads
\begin{align}
\label{Eq:RieszPair}
\Fourier\big\{\|\cdot\|^{\lambda}\big\}(\bw)
=\frac{2^{\lambda+d}\,\pi^{d/2}\,\Gamma\big(\tfrac{\lambda+d}{2}\big)}
{\Gamma\big(-\tfrac{\lambda}{2}\big)}\,\|\bw\|^{-\lambda-d},
\end{align}
which holds by analytic continuation for $\lambda \notin \{0,2,4,\dots\}\cup\{-d,-d-2,\dots\}$, where the excluded values are those for which the pair degenerates and logarithmic factors appear \big(see \cite[Chapter~II]{Gelfand-Shilov1964} and \cite[Section~2.4.3]{Grafakos2008}\big). The choices $\lambda=\gamma-d$
and $\lambda=2\gamma-d=2H$ yield
\begin{align}
\label{Eq:fBfGreen}
\rho_\Lop(\V x)&=C_{\gamma,d}\,\|\V x\|^{\gamma-d}
\qquad \mbox{ with } \qquad
C_{\gamma,d}=\frac{\Gamma\big(\tfrac{d-\gamma}{2}\big)}
{2^{\gamma}\,\pi^{d/2}\,\Gamma\big(\tfrac{\gamma}{2}\big)}\\
\label{Eq:fBfGreenB}
\rho_{\Op B}(\V x)&=c_{\gamma,d}\,\|\V x\|^{2H}
\qquad \mbox{ with } \qquad
c_{\gamma,d}=\frac{\Gamma(-H)}{4^{\gamma}\,\pi^{d/2}\,\Gamma(\gamma)},
\end{align}
where $\rho_\Lop=\Fourier^{-1}\{\|\cdot\|^{-\gamma}\}$ and
$\rho_{\Op B}=\Fourier^{-1}\{\|\cdot\|^{-2\gamma}\}$ are the canonical (isotropic) Green's
functions of $\Lop$ and $\Op B$. The formula for $\rho_{\Op B}$ is valid throughout
\eqref{Eq:fBfRange}, since $2H$ avoids the even integers exactly when $H \notin \N$. The
formula for $\rho_\Lop$ additionally requires $\gamma-d \notin \{0,2,4,\dots\}$; at those
values, $\rho_\Lop$ acquires a logarithmic factor, which affects the innovation kernel
below but neither $\rho_{\Op B}$ nor the covariance.

Both Green's functions are powers of the Euclidean norm, and the boundary corrections that
follow involve their Taylor polynomials at the origin. We therefore record the closed
form of the directional derivatives of $f_\lambda=\|\cdot\|^{\lambda}$,
\begin{align}
\label{Eq:DirDeriv}
(\V h\cdot\nabla)^m f_\lambda(\V u)
=m!\sum_{l=0}^{\lfloor m/2\rfloor}\binom{\lambda/2}{m-l}\binom{m-l}{l}
\big(2\langle \V u,\V h\rangle\big)^{m-2l}\|\V h\|^{2l}\|\V u\|^{\lambda-2(m-l)},
\end{align}
where $\binom{\lambda/2}{j}$ is the generalized binomial coefficient. Indeed,
$\|\V u+t\V h\|^{2}=\|\V u\|^2+2t\langle \V u,\V h\rangle+t^2\|\V h\|^2$, so that
$t \mapsto f_\lambda(\V u+t\V h)$ is a power of a quadratic polynomial in $t$; the binomial
series $(1+v)^{\lambda/2}=\sum_{j\ge0}\binom{\lambda/2}{j}v^j$, applied to
$v=\big(2t\langle \V u,\V h\rangle+t^2\|\V h\|^2\big)/\|\V u\|^2$ and followed by the
expansion of $\big(2\langle \V u,\V h\rangle+t\|\V h\|^2\big)^j$, gives \eqref{Eq:DirDeriv}
with $l=m-j$ counting the powers of $t$ that come from the quadratic term. An equivalent
form is
\begin{align*}
(\V h\cdot\nabla)^m f_\lambda(\V u)
=(-1)^m\,m!\,\|\V h\|^m\,\|\V u\|^{\lambda-m}\,C^{(-\lambda/2)}_m(\cos\theta),
\end{align*}
where $\cos\theta=\langle \V u,\V h\rangle/(\|\V u\|\,\|\V h\|)$ and $C^{(\alpha)}_m$ is the
Gegenbauer polynomial of degree $m$.

The three objects of Theorem \ref{Theo:ReproKernel} retain a compact form at every order.
Since the Riesz kernels are radial, we have that
$q_{\V k}(\V y)=\langle \phi_{\V k},\rho_\Lop(\cdot-\V y)\rangle
=(\partial^{\V k}\rho_\Lop)(-\V y)$, so that Item \ref{Item:KernelInverse} identifies the
impulse response of the innovation model as the Taylor remainder
\begin{align}
\label{Eq:fBfKernel}
g_{\V \phi}(\V x,\V y)=\rho_\Lop(\V x-\V y)
-\sum_{|\V k|\le k_\gamma}\frac{\V x^{\V k}}{\V k!}\,(\partial^{\V k}\rho_\Lop)(-\V y),
\end{align}
in which the subtraction of the Taylor polynomial of $\V x \mapsto \rho_\Lop(\V x-\V y)$ at
the origin is the action of the projector $\Op P^\ast_{\V \phi}$, exactly as in the
one-dimensional case of Appendix \ref{App:fBm}. The multinomial theorem
$\sum_{|\V k|=m}\frac{\V x^{\V k}}{\V k!}\partial^{\V k}=\frac{1}{m!}(\V x\cdot\nabla)^m$,
together with \eqref{Eq:DirDeriv} for $\lambda=\gamma-d$, turns \eqref{Eq:fBfKernel} into
the closed form
\begin{align}
\label{Eq:fBfKernelExplicit}
g_{\V \phi}(\V x,\V y)=C_{\gamma,d}\Bigg(&\|\V x-\V y\|^{\gamma-d}\nonumber\\
&-\sum_{m=0}^{k_\gamma}\sum_{l=0}^{\lfloor m/2\rfloor}
\binom{\tfrac{\gamma-d}{2}}{m-l}\binom{m-l}{l}
\big(\!-2\langle \V x,\V y\rangle\big)^{m-2l}\nonumber\\
&\qquad \qquad \qquad \qquad \times\|\V x\|^{2l}\,\|\V y\|^{\gamma-d-2(m-l)}\Bigg),
\end{align}
in which the factorials of \eqref{Eq:DirDeriv} cancel those of the Taylor expansion. The
fractional Brownian field of order $\gamma$ is the corresponding
process $\gRand B_\gamma=\Op L^{-1}_{\V \phi}\{\gRand W\}$; that is, the solution of the
stochastic partial differential equation
\begin{align}
\label{Eq:fBfSDE}
(-\Delta)^{\gamma/2}\{\gRand B_\gamma\}=\gRand W
\qquad \mbox{ subject to } \qquad \partial^{\V k}\gRand B_\gamma(\V 0)=0,
\quad |\V k|\le k_\gamma.
\end{align}
Equivalently, it is the stochastic integral
\begin{align}
\label{Eq:fBfMandelbrot}
\gRand B_\gamma(\V x)=\int_{\R^d}g_{\V \phi}(\V x,\V z)\,\gRand W(\V z)\,\dint \V z
\end{align}
with the kernel \eqref{Eq:fBfKernelExplicit}, which is the multidimensional counterpart, at
arbitrary order, of Mandelbrot's representation \eqref{Eq:fBMH}. For $k_\gamma=0$, the
correction reduces to the single term $\|\V z\|^{\gamma-d}$, so that
$\gRand B_\gamma(\V x)=C_{\gamma,d}\int_{\R^d}\big(\|\V x-\V z\|^{\gamma-d}
-\|\V z\|^{\gamma-d}\big)\gRand W(\V z)\,\dint \V z$, which is L\'evy's field. Let us note
that the Riesz kernel is symmetric, so that the case $d=1$ returns the well-balanced version
of fBm rather than the causal one of \eqref{Eq:fBMH}; the two differ by a normalization
factor only, since fBm is determined by $H$ up to a multiplicative constant.
Likewise, $v_{\V k}=\Op A\{\phi_{\V k}\}=(-1)^{|\V k|}\,\partial^{\V k}\rho_{\Op B}$ and
$r_{\V k,\V l}=\pm(\partial^{\V k+\V l}\rho_{\Op B})(\V 0)=0$, because
$\partial^{\V k+\V l}\rho_{\Op B}$ is homogeneous of degree
$2H-|\V k|-|\V l|\ge 2H-2k_\gamma>0$ and therefore vanishes at the origin.

The conditioning has a transparent Fourier signature, and this is where the multi-index sums
telescope. In the Fourier domain, the projector $\Op P_{\V \phi}$ subtracts from $\hat \varphi$ its
Taylor polynomial of degree $k_\gamma$ at $\bw=\V 0$. For the sampling
functional $\delta_{\V x}$, whose Fourier transform is $\ee^{-\jj \langle \bw,\V x\rangle}$,
the multinomial theorem collapses that polynomial---a priori, a sum over the $N_0$
multi-indices---into a function of the single scalar $\langle \bw,\V x\rangle$, so that
$\Fourier\{\Op P_{\V \phi}\delta_{\V x}\}(\bw)
=E_{k_\gamma}\big(-\jj \langle \bw,\V x\rangle\big)$ with the truncated exponential
\begin{align}
\label{Eq:TruncExp}
E_K(t)\eqdef \ee^{t}-\sum_{m=0}^{K}\frac{t^m}{m!}.
\end{align}
By Parseval's relation, the identity
$a_{\V \phi}(\V x,\V y)=\langle \Op A\,\Op P_{\V \phi}\{\delta_{\V y}\},
\Op P_{\V \phi}\{\delta_{\V x}\}\rangle$ then yields the Fourier form of the covariance
function
\begin{align}
\label{Eq:fBfFourier}
a_{\V \phi}(\V x,\V y)=\frac{1}{(2\pi)^d}\int_{\R^d}
\frac{E_{k_\gamma}\big(-\jj \langle \bw,\V x\rangle\big)\,
\overline{E_{k_\gamma}\big(-\jj \langle \bw,\V y\rangle\big)}}{\|\bw\|^{2\gamma}}\,\dint \bw,
\end{align}
which is the specialization to sampling functionals of the quadratic form
\eqref{Eq:fBmQuadratic} of Appendix \ref{App:fBm}, and from which the qualitative
properties of the field can be read off directly.

\begin{proposition}[Covariance of the fractional Brownian fields]
\label{Theo:fBfCov}
Under the condition \eqref{Eq:fBfRange}, the covariance function of the field
\eqref{Eq:fBfSDE} is the double Taylor remainder
\begin{align}
\label{Eq:fBfCov}
a_{\V \phi}(\V x,\V y)=c_{\gamma,d}\Bigg(&\|\V x-\V y\|^{2H}\nonumber\\
&-\sum_{m=0}^{k_\gamma}\frac{(-1)^m}{m!}
\Big[(\V x\cdot\nabla)^m f_{2H}(\V y)+(\V y\cdot\nabla)^m f_{2H}(\V x)\Big]\Bigg),
\end{align}
where $f_{2H}=\|\cdot\|^{2H}$, where $(\V h\cdot\nabla)^m$ denotes the $m$-fold directional
derivative along $\V h$, and with $c_{\gamma,d}$ given by \eqref{Eq:fBfGreenB}. These
derivatives are given in closed form by \eqref{Eq:DirDeriv} with $\lambda=2H$, and the
factorials cancel upon substitution, so that \eqref{Eq:fBfCov} is an elementary finite
combination of $\|\V x\|$, $\|\V y\|$, $\langle \V x,\V y\rangle$ and
$\|\V x-\V y\|^{2H}$. The field is isotropic and self-similar of exponent $H$ at every
order, with
\begin{align}
\label{Eq:fBfSigma}
\Exp\big\{|\gRand B_\gamma(\V x)|^2\big\}=\sigma^2_{\gamma,d}\,\|\V x\|^{2H},
\qquad
\sigma^2_{\gamma,d}=-2\,(-1)^{k_\gamma}\binom{2H-1}{k_\gamma}\,c_{\gamma,d}>0.
\end{align}
For $k_\gamma=0$, \eqref{Eq:fBfCov} reduces to the classical form
\begin{align}
\label{Eq:fBfCovLevy}
a_{\V \phi}(\V x,\V y)
=\frac{\sigma^2_{\gamma,d}}{2}\Big(\|\V x\|^{2H}+\|\V y\|^{2H}-\|\V x-\V y\|^{2H}\Big),
\qquad
\sigma^2_{\gamma,d}=\frac{2\,\Gamma(1-H)}{H\,4^{\gamma}\,\pi^{d/2}\,\Gamma(\gamma)}.
\end{align}
\end{proposition}
\begin{proof}
Formula \eqref{Eq:fBfCov} is \eqref{Eq:kernelrepro} of Theorem \ref{Theo:ReproKernel} with
$r_{\V k,\V l}=0$ and
$p_{\V k}(\V x)\,v_{\V k}(\V y)=c_{\gamma,d}\frac{(-\V x)^{\V k}}{\V k!}\,\partial^{\V k}f_{2H}(\V y)$,
in which the multinomial theorem,
\begin{align*}
\sum_{|\V k|=m}\frac{\V h^{\V k}}{\V k!}\,\partial^{\V k}
=\frac{1}{m!}\,(\V h\cdot\nabla)^m,
\end{align*}
collapses the sums over the $N_0$ multi-indices into the $(k_\gamma+1)$ directional terms.
The Fourier form \eqref{Eq:fBfFourier} then delivers the qualitative properties: the
substitution $\bw \mapsto \M R\bw$ with $\M R$ a rotation gives the isotropy; the dilation
$\bw \mapsto \bw/s$ gives $a_{\V \phi}(s\V x,s\V y)=s^{2H}\,a_{\V \phi}(\V x,\V y)$
(self-similarity of exponent $H$); and the bounds $E_{k_\gamma}(t)=O(t^{k_\gamma+1})$ at the
origin and $E_{k_\gamma}(t)=O(t^{k_\gamma})$ at infinity (for imaginary $t$) show that the
integral is finite exactly when $k_\gamma<H<k_\gamma+1$, which is the range
\eqref{Eq:fBfRange}. The finiteness of the diagonal makes the pointwise reading legitimate,
by Item \ref{Item:RKHS} of Theorem \ref{Theo:ReproKernel}.

It remains to evaluate the covariance. Let us set $\V y=\V x$ in \eqref{Eq:fBfCov}: the first
term vanishes, while the homogeneity of $f_{2H}$ gives, along the ray,
\begin{align*}
(\V x\cdot\nabla)^m f_{2H}(\V x)
&=\frac{\dint^m}{\dint s^m}\Big[(1+s)^{2H}\Big]_{s=0}\,f_{2H}(\V x)\\
&=2H(2H-1)\cdots(2H-m+1)\,f_{2H}(\V x),
\end{align*}
which is also the case $\V h=\V x$ of \eqref{Eq:DirDeriv}, so that
\begin{align*}
a_{\V \phi}(\V x,\V x)
=-2\,c_{\gamma,d}\sum_{m=0}^{k_\gamma}(-1)^m\binom{2H}{m}\,\|\V x\|^{2H}.
\end{align*}
The alternating-sum identity
$\sum_{m=0}^{K}(-1)^m\binom{\alpha}{m}=(-1)^{K}\binom{\alpha-1}{K}$ then gives
\eqref{Eq:fBfSigma}. The result is positive: the factors $2H-1,\dots,2H-k_\gamma$ of the
generalized binomial coefficient are positive because $2H>2k_\gamma$, while the sign of
$\Gamma(-H)$, and hence of $c_{\gamma,d}$, is $(-1)^{k_\gamma+1}$.
\end{proof}
For $d=1$, Proposition \ref{Theo:fBfCov} reduces to Proposition \ref{Theo:fBmCov} of
Appendix \ref{App:fBm}, whose variance formula \eqref{Eq:fBmVariance} the alternating-sum
identity above condenses into the single generalized binomial coefficient of
\eqref{Eq:fBfSigma}, as one verifies with the duplication and reflection formulas of the
Gamma function. In particular, the $d=1$ case of \eqref{Eq:fBfCov} is the covariance of
the $n$th-order fractional Brownian motion of Perrin et al.\ \cite{Perrin2001} with
$n=k_\gamma+1$. By contrast, for $d\ge 2$ and $H>1$, the classical literature treats the
isotropic self-similar Gaussian fields as generalized random fields that are defined only
modulo the polynomials of degree $k_\gamma$, whether through the spectral representation
\cite{Gelfand-Vilenkin1964,Dobrushin1979,Lodhia2016} or through the generalized
covariances of Matheron's intrinsic random functions of order $k$ \cite{Matheron1973}.
The conditioning by $\V \phi$ removes this ambiguity, and this is what gives access to
the pointwise formula \eqref{Eq:fBfCov}. To the best of our knowledge, the explicit
description that \eqref{Eq:fBfKernelExplicit}, \eqref{Eq:fBfMandelbrot} and
\eqref{Eq:fBfCov} provide for $d \ge 2$ and $k_\gamma \ge 1$ is new. The operator-based characterizations of the fractional Brownian
fields \cite{Tafti2010} specify these processes through their whitening operator and their
spectral measure rather than through a pointwise covariance, and the closed forms that are
available in the literature are the case $k_\gamma=0$ of \eqref{Eq:fBfCovLevy}, which is
L\'evy's field, and the one-dimensional formula of \cite{Perrin2001}.

Anticipating Section \ref{Sec:MMSESplines}, Theorem \ref{Theo:MMSE} will identify the MMSE
estimator of the field from the noisy samples $y_m=\gRand B_\gamma(\V x_m)+\epsilon_m$ as the constrained spline
\begin{align}
\label{Eq:fBfSpline}
\tilde g_\lambda=\sum_{m=1}^M a_m\, a_{\V \phi}(\cdot,\V x_m),
\qquad \M a=(\widetilde{\M G}+\lambda \M I)^{-1}\M y, \qquad \lambda=\sigma^2,
\end{align}
with $[\widetilde{\M G}]_{m,n}=a_{\V \phi}(\V x_m,\V x_n)$. These are the polyharmonic (or
fractional thin-plate) splines of Duchon \cite{Duchon1977}, whose radial basis functions
$\|\cdot-\V x_m\|^{2H}$ are $\Pi_{k_\gamma}$-conditionally positive-definite in the sense of
Definition \ref{Def:CondPositiveDef1}, the boundary corrections being carried by
$a_{\V \phi}$. The excluded integer values of $H$ in \eqref{Eq:fBfRange} are precisely those
of the classical thin-plate splines, whose radial functions (e.g., $\|\V x\|^{2}\log\|\V x\|$
for $\gamma=2$ and $d=2$) acquire the logarithmic factor. The interpolation counterpart
(Corollary \ref{Theo:Unser}) applies as soon as the deterministic
observations $\partial^{\V k}\gRand B_\gamma(\V 0)=0$, $|\V k|\le k_\gamma$, are part of the
data.

\section{Generalized Splines as MMSE Estimators of GGPs}
\label{Sec:MMSESplines}
We now have all the elements to
make the connection between the minimum-error estimation of a signal under the Gaussian hypothesis---also known as the {\em Wiener estimator}---and the spline reconstruction techniques of \cite[Section~\ref{II-Sec:InverseProblems}]{Unser2026b}.
To that end, we consider the linear measurement model
\begin{align}
\label{Eq:Wienerproblem}
g \mapsto \M y=\V \nu(g) + \V \epsilon \in \R^M,
\end{align}
where $g$ (our signal) is a realization of a Gaussian process $\gRand g \sim \Spc N(0,\Op A)$ in $\Spc S'$  
 and where $\V \nu=(\nu_m)$ with $\nu_1,\dots, \nu_M\in \Spc H'_\Lop$ is a continuous linear measurement operator $\V \nu: \Spc H_\Lop \to \R^M$.
The second component of the model, $\V \epsilon \in \R^M$, is an additive disturbance term (discrete measurement noise) that does not depend on the signal and whose components are i.i.d.\ Gaussian with zero mean and variance
$\sigma^2$.

\subsection{Estimation of Non-Degenerate GGPs}
\label{Sec:EstimConventional}
In the interest of clarity, we start with the ``conventional'' scenario $\Op A=(\Lop^\ast\Lop)^{-1}: \Spc S \to \Spc S'$
with $\Lop$ invertible ($N_0=0$), for which we have explicitly solved
the deterministic inverse problem in \cite[Section~\ref{II-Sec:Tikhonov}]{Unser2026b}. Specifically, we know that the {\em Tikhonov estimator} of $g$ in \eqref{Eq:Wienerproblem} from the noisy measurements $\M y$ is the generalized smoothing spline
\begin{align}
g_\lambda&=\arg \min_{g \in \Spc H_\Lop}\left( \|\M y - \V \nu(g)\|^2_2 + \lambda \|\Lop g\|^2_{L_2}\right) \nonumber\\
&=\sum_{m=1}^M a_m \Op A\{\nu_m\}\quad \mbox{ with } \quad (a_m)=(\M G + \lambda \M I)^{-1} \M y, 
\label{Eq:GSmoothingSpline}
\end{align}
where the matrix $\M G \in \R^{M \times M}$ with $[\M G]_{m,n}=\langle \Op A \nu_m, \nu_n\rangle=\langle \nu_m,\nu_n\rangle_{\Spc H'_\Lop}$ is symmetric positive (semi-)definite.

Let us note that the matrix inverse in \eqref{Eq:GSmoothingSpline} is always well-defined for $\lambda>0$. 
For the limit case of an exact fit with $\lambda=0$, we need the additional assumption that $\M G$ is positive-definite, which is equivalent to the $\nu_m$ being linearly independent. When this latter assumption fails, it suffices to reduce the (noise-free) data set by selecting a subset of $r={\rm rank}(\M G)<M$ measurements $\widetilde{\M y}\in \R^r$ such that the reduced Gram matrix $\widetilde{\M G}\in \R^{r \times r}$ is of full rank. This can be implemented systematically by replacing the ill-posed inverse of $\M G$ by the generalized Moore-Penrose inverse $\M G^+$ (see \cite[pp. 284-295]{BenIsrael2003}). Since
$\M A^+=\M A^{-1}$ when $\M A$ is invertible, we can therefore replace the right-hand-side of \eqref{Eq:GSmoothingSpline} by its foolproof variant $\M a=(\M G + \lambda \M I)^{+} \M y$.

We now turn to the Bayesian counterpart: the MMSE estimator of $g(\varphi)$, for a fixed
$\varphi \in \Spc H'_\Lop$, given the measurements $\M y$. In probabilistic terms, the model
\eqref{Eq:Wienerproblem} reads $\VRand Y=\VRand Y_0+\V \epsilon$ with
$\VRand Y_0=\big(\gRand g(\nu_1),\dots,\gRand g(\nu_M)\big)$ and
$\V \epsilon \sim \Spc N(\V 0,\sigma^2 \M I)$ independent of $\gRand g$, while the quantity
to estimate is the scalar random variable $\Rand X=\gRand g(\varphi)$. By
Theorem \ref{Theo:GGPcondi}, which covers the test functions
$\varphi,\nu_1,\dots,\nu_M \in \Spc H'_\Lop$, the vector $(\Rand X,\VRand Y_0)$ is Gaussian
with the covariances $\langle \Op A\varphi,\varphi\rangle$,
$\langle \Op A\varphi,\nu_m\rangle$, and $\langle \Op A\nu_m,\nu_n\rangle$; the independence
of the noise then gives
\begin{align}
\label{Eq:AugmentedCovarZ}
\VRand z=(\Rand X,\VRand Y) \sim \Spc N(\V 0,\M C_{\VRand z}),
\qquad
\M C_{\VRand z}=\left(\begin{array}{cc}
\sigma_\Rand x^2  &   \M c^\Tr    \\[0.5ex]
\M c  &   \M G + \sigma^2 \M I 
\end{array}
\right),
\end{align}
with $\sigma^2_{\Rand x}=\langle \Op A\varphi,\varphi\rangle$,
$[\M c]_m=\langle \Op A\{\nu_m\},\varphi\rangle$, and
$[\M G]_{m,n}=\langle \Op A\nu_m,\nu_n\rangle$. Two features of \eqref{Eq:AugmentedCovarZ}
are noteworthy: the vector $\M c=\V \nu^\ast(\varphi)$ collects the evaluations at $\varphi$
of the Riesz representers $\V \nu^\ast=\big(\Op A\{\nu_1\},\dots,\Op A\{\nu_M\}\big)$ of the
measurement functionals, and $\M G$ is the very Gram matrix of \eqref{Eq:GSmoothingSpline}.
The conditional law of one block of a Gaussian vector given the other block is Gaussian, with
mean and variance that are read off the partitioned covariance
\cite[Eqs.~(2.81)--(2.82)]{Bishop2006}. Here, $p_{\Rand X|\VRand Y=\M y}$ is the univariate
Gaussian with mean
\begin{align}
\Exp\{\Rand X|\M y\}=\M c^\Tr (\M G+\sigma^2 \M I)^{-1}\M y
=\V \nu^\ast(\varphi)^\Tr (\M G + \sigma^2 \M I)^{-1} \M y
\end{align}
and variance
$\sigma^2_{\Rand X|\M y}=\langle \Op A \varphi,\varphi\rangle
-\V \nu^\ast(\varphi)^\Tr (\M G + \sigma^2 \M I)^{-1}\V \nu^\ast(\varphi)$. Since the
conditional mean is the MMSE estimator whatever the distribution, the comparison with
\eqref{Eq:GSmoothingSpline} yields the following.
\begin{theorem}[Generalized Tikhonov-Wahba-Wiener theorem] 
\label{Theo:GenGaussMarkov}
We consider the following setting.
\begin{itemize}
\item The signal $g=\{g(\varphi)\}_{\varphi \in \Spc S} \in \Spc S'$ to be recovered  is a realization of a generalized Gaussian process $\gRand g\sim \Spc N(0,\Op C_\Lop) $ in $\Spc S'$ with $\Op C_\Lop$ positive definite.
\item There exists an invertible whitening/regularization operator $\Lop$ such that $\Op C_\Lop=\Op A=(\Lop^\ast \Lop)^{-1}: \Spc S \to \Spc S'$. This operator induces the native Hilbert space $\Spc H_\Lop=\overline{(\Spc S,\|\Lop\, \cdot\|_{L_2})}$ of Section \ref{Sec:NativeRecap}, whose continuous dual $\Spc H'_\Lop$ is equipped with the inner product $\langle \nu_1,\nu_2\rangle_{\Spc H'_\Lop}=\langle \Op A \nu_1,\nu_2\rangle$.
\item The observed data are $\M y=\V \nu(g) + \V \epsilon \in \R^M$,
with $\V \nu=(\nu_1, \dots,\nu_M)\in (\Spc H'_\Lop)^M$ and $\V \epsilon$ a realization of a white Gaussian noise  with variance $\sigma^2$.
\end{itemize}
Then, for any $\varphi\in \Spc H'_\Lop\supset \Spc S$, the MMSE estimator of
$g(\varphi)$ given $\M y$
is 
\begin{align}
\label{Eq:MMSEsol0}
g_{{\rm MMSE}}(\varphi|\M y)=\Exp\{\gRand g( \varphi)|\M y\}%
&=g_{\lambda}( \varphi)\quad \mbox{ with }\quad  \lambda=\sigma^2,
\end{align}
where
$g_{\lambda}\in \Spc H_\Lop$ is the generalized smoothing spline defined by \eqref{Eq:GSmoothingSpline}.

\end{theorem}

Since the result in Theorem \ref{Theo:GenGaussMarkov} holds for any $\varphi\in  \Spc S$, it gives a precise meaning to the statement that the generalized spline $g_\lambda=\{g_\lambda(\varphi): \varphi \in \Spc S\} \in  \Spc S'$ with $\lambda=\sigma^2$ is the optimal estimate of the generalized Gaussian process
$\gRand G=\{\gRand g(\varphi): \varphi \in \Spc S\}$ in $\Spc S'$ given the measurements $\M y$. As for the case where $\Spc H_\Lop$ is a RKHS of ordinary functions $f: \R^d \to \R$, we have that $\delta_{\V x_0} \in \Spc H'_\Lop$ for any 
$\V x_0 \in \R^d$. Consequently, we can plug $\varphi=\delta_{\V x_0}$ in \eqref{Eq:MMSEsol0}, which yields Wahba's equivalence for classic  Gaussian processes $\gRand G=\{\gRand g(\V x)\}_{\V x \in \R^d}$ whitened by $\Lop$; namely, that $\Exp\{\gRand g(\V x_0)|\M y\}=g_{\lambda}(\V x_0)$
with $\lambda=\sigma^2$.

\subsection{Estimation of GGPs with Non-Trivial Null Space}

The variational side of the estimation is governed by the representer theorem
\cite[Theorem~\ref{II-Theo:GeneralRepSemiHilbert}]{Unser2026b}, which states that the solution of a
regularized inverse problem over $\Spc H_\Lop$ with measurement functionals
$\nu_1,\dots,\nu_M \in \Spc H'_\Lop$ lies in the finite-dimensional space that is spanned
by $\{\Op A_{\V \phi}\{\nu_m\}\}_{m=1}^M \cup \{p_n\}_{n=1}^{N_0}$. This requires the
measurement functionals $\V \nu=(\nu_1,\dots,\nu_M)$ to be {\em $\Lop$-compatible}
\cite[Definition~\ref{II-Def:AdmisPoints}]{Unser2026b}: the $\nu_m$ are linearly independent
elements of $\Spc H'_\Lop$, and the matrix $\M P \in \R^{M \times N_0}$ with
$[\M P]_{m,n}=\langle \nu_m,p_n\rangle$ is of rank $N_0$.

We now present two variants of Theorem \ref{Theo:GenGaussMarkov} for GGPs that are whitened by a $\V p$-admissible operator $\Lop$.
The first is the direct generalization of a weaker variant of \eqref{Eq:MMSEsol0} for ordinary stochastic processes that states that $\Lop$-smoothing splines yield
the {\em best linear unbiased} (BLU) estimator of a stochastic process $\{\gRand X(t)\}$ whitened by $\Lop$ under the assumption of a fixed-effect model that includes an additional deterministic drift $p(\cdot) \in \Spc N_\V p$ \cite{Wahba1990,Berlinet2004}. In our formalism, the generalized version of this model reads
\begin{align}
\label{Eq:FixedEffect}
\begin{cases}\gRand G&=\sum_{n=1}^{N_0}  \theta_n p_n + \gRand X\\
\Rand Y_m&=\gRand G(\nu_m) + \epsilon_m,\quad m=1,\dots,M,
\end{cases}
\end{align}
where the $\theta_n$ are fixed but unknown, $\V \epsilon=(\epsilon_m)\sim \Spc N(0,\sigma^2\M I)$, and where $\Rand X\sim \Spc N(0,\Op A_{\V \phi})$ is a GGP
whose covariance operator $\Op A_{\V \phi}: \Spc S \to \Spc S'$ is the Gram operator \eqref{III-Eq:Aphi} that is built from $\Op A=(\Lop^\ast \Lop)^{-1}$ where $\Lop$ is ${\V p}$-admissible with $\V p=(p_1, \dots,p_{N_0})$.

The BLU estimator of $\gRand G(\varphi)$ in \eqref {Eq:FixedEffect} given the observed data $\VRand Y=\M y$ has the form
$g_{\rm BLU}(\varphi|\M y)=\V \beta^\Tr \M y$ (linearity), with its coefficient-vector $\V \beta=\V \beta(\varphi) \in \R^M$ being selected such as to minimize $\Exp\{\big|g_{\rm BLU}(\varphi|\M y)- \gRand G(\varphi)\big|^2\}$, subject to the zero-bias constraint $\Exp\{g_{\rm BLU}(\varphi|\M y)- \gRand G(\varphi)|\V \theta\}=0$ for every $\V \theta \in \R^{N_0}$. The relation of this estimator to splines is stated in Theorem \ref{Theo:WahbaKimeldorf}. In the classical setting, the equivalence is stated in \cite[Theorems 1.5.1--1.5.2]{Wahba1990}, with the proof deferred to the explicit computation of \cite{Kimeldorf1971}; an alternative route derives the smoothing spline as the limit of Bayes estimators under an improper prior on the null-space coefficients \cite[Theorem 1.5.3]{Wahba1990}. The proof below is self-contained: it reduces the equivalence to the symmetry of a single matrix, without any explicit inversion.
\begin{theorem}[Kimeldorf--Wahba \cite{Kimeldorf1971}]
\label{Theo:WahbaKimeldorf}
Let $\V \nu$ be $\Lop$-compatible in the sense defined above, let $\sigma^2>0$, and let $g_\lambda$ be the generalized smoothing spline defined by
\begin{align}
g_\lambda=\arg \min_{g \in \Spc H_\Lop}\left( \|\M y - \V \nu(g)\|^2_2 + \lambda \|\Lop g\|^2_{L_2}\right).
\end{align}
Then, for any $\varphi \in \Spc H'_\Lop$ and $\lambda=\sigma^2$, $g_\lambda(\varphi)=g_{\rm BLU}(\varphi|\M y)$ where $g_{\rm BLU}(\varphi|\M y)$ is the {\rm best linear unbiased estimator} of $\gRand G(\varphi)$ in \eqref {Eq:FixedEffect} given $\M y$.
\end{theorem}

\begin{proof}
Throughout, we use the matrix $\M P=\big[\langle \nu_m,p_n\rangle\big] \in \R^{M \times N_0}$
(of rank $N_0$, by $\Lop$-compatibility), the Gram matrix
$\widetilde{\M G}=\big[\langle \Op A_{\V \phi}\{\nu_m\},\nu_n\rangle\big] \in \R^{M \times M}$,
and the two $\varphi$-dependent vectors
$\widetilde{\M c}=\big(\langle \Op A_{\V \phi}\{\nu_m\},\varphi\rangle\big) \in \R^{M}$ and
$\V p(\varphi)=\big(\langle p_n,\varphi\rangle\big) \in \R^{N_0}$.

\medskip
\emph{Normal equations of the BLU estimator}:
In the model \eqref{Eq:FixedEffect}, we have that
$\VRand Y=\M P\V \theta+\V \nu(\Rand X)+\V \epsilon$ and
$\Exp\{\gRand G(\varphi)\}=\V \theta^\Tr \V p(\varphi)$, so that the bias of
$g_{\rm BLU}(\varphi|\M y)=\V \beta^\Tr \M y$ is
$\big(\M P^\Tr \V \beta - \V p(\varphi)\big)^\Tr \V \theta$. The requirement that it
vanishes for every $\V \theta \in \R^{N_0}$ is the linear constraint
$\M P^\Tr \V \beta=\V p(\varphi)$, which is feasible because $\M P$ has rank $N_0$. Under
that constraint, the deterministic part cancels out from the error, which reduces to
$\V \beta^\Tr\big(\V \nu(\Rand X)+\V \epsilon\big)-\Rand X(\varphi)$, with variance
\begin{align*}
\V \beta^\Tr \big(\widetilde{\M G}+\sigma^2 \M I\big)\V \beta
-2\,\V \beta^\Tr \widetilde{\M c}+\langle \Op A_{\V \phi}\varphi,\varphi\rangle,
\end{align*}
where the second-order statistics of $\Rand X$ are given by \eqref{Eq:Cov2}. The
minimization of this quadratic form subject to $\M P^\Tr \V \beta=\V p(\varphi)$ yields, by
Lagrange duality, the linear system
\begin{align}
\label{Eq:BLUSystem}
\M K \begin{pmatrix} \V \beta \\ \V \gamma \end{pmatrix}
=\begin{pmatrix} \widetilde{\M c} \\ \V p(\varphi) \end{pmatrix},
\qquad
\M K \eqdef \begin{pmatrix} \widetilde{\M G}+\sigma^2 \M I & \M P \\ \M P^\Tr & \M 0
\end{pmatrix},
\end{align}
whose symmetric matrix $\M K$ is invertible because $\widetilde{\M G}+\sigma^2\M I$ is
positive-definite and $\M P$ is of rank $N_0$.

\medskip
\emph{Normal equations of the smoothing spline}:
By the representer theorem of \cite{Unser2026b}, as recalled in Section \ref{Sec:NativeRecap}, the minimizer $g_\lambda$ is
unique and can be written as
$g_\lambda=\sum_{m=1}^M a_m \Op A_{\V \phi}\{\nu_m\}+\sum_{n=1}^{N_0}b_n p_n$.
The two building
blocks of the loss are then
$\V \nu(g_\lambda)=\widetilde{\M G}\M a+\M P \M b$ and
$\|\Lop g_\lambda\|^2_{L_2}=\M a^\Tr \widetilde{\M G}\M a$, where the second identity
follows from $\Lop \Op A_{\V \phi}=\Op L^{-1\ast}_{\V \phi}$ and the factorization
\eqref{III-Eq:Factorization}. The zeroing of the gradient of
$\|\M y-\widetilde{\M G}\M a-\M P\M b\|^2_2+\lambda \M a^\Tr\widetilde{\M G}\M a$ gives the
stationarity conditions
$\widetilde{\M G}\big[(\widetilde{\M G}+\lambda \M I)\M a+\M P\M b-\M y\big]=\V 0$ and
$\M P^\Tr\big[\widetilde{\M G}\M a+\M P\M b-\M y\big]=\V 0$, and one verifies directly that
these are satisfied by the solution of
\begin{align}
\label{Eq:SplineSystem}
\M K \begin{pmatrix} \M a \\ \M b \end{pmatrix}
=\begin{pmatrix} \M y \\ \V 0 \end{pmatrix}
\qquad \mbox{with } \lambda=\sigma^2,
\end{align}
since the latter gives $(\widetilde{\M G}+\lambda\M I)\M a+\M P\M b=\M y$ and
$\M P^\Tr \M a=\V 0$, hence
$\M P^\Tr\big[\widetilde{\M G}\M a+\M P\M b-\M y\big]=-\lambda \M P^\Tr\M a=\V 0$. The loss
is convex, so that stationarity suffices, and the evaluation of the solution reads
$g_\lambda(\varphi)=\M a^\Tr \widetilde{\M c}+\M b^\Tr \V p(\varphi)$.

\medskip
\emph{Equivalence}:
The two estimators are driven by the same symmetric invertible matrix $\M K$, with the
roles of the data and of the test function exchanged. Specifically,
\begin{align*}
g_\lambda(\varphi)
=\begin{pmatrix} \widetilde{\M c} \\ \V p(\varphi)\end{pmatrix}^{\!\Tr}
\M K^{-1}\begin{pmatrix} \M y \\ \V 0\end{pmatrix}
=\begin{pmatrix} \M y \\ \V 0\end{pmatrix}^{\!\Tr}
\M K^{-1}\begin{pmatrix} \widetilde{\M c} \\ \V p(\varphi)\end{pmatrix}
=\V \beta^\Tr \M y=g_{\rm BLU}(\varphi|\M y),
\end{align*}
where the middle equality is the symmetry of $\M K^{-1}$.
\end{proof}

Theorem \ref{Theo:WahbaKimeldorf} is remarkable in that it shows the
complete equivalence between generalized smoothing splines and BLU estimators for GGPs.
The same theorem clarifies the status of the polyharmonic splines of Duchon \cite{Duchon1977} for $d\ge3$, whose polynomial part ranges over $\Pi_{m-1}\supsetneq\Spc N_{\V p}$ \cite[Remark~\ref{I-Rem:Duchon}]{Unser2026a}: the components in the excess polynomials are fixed effects in the sense of \eqref{Eq:FixedEffect}, so that these splines are BLU estimators, but they are not the MMSE estimator of any GGP; conditioning the process on $\Pi_{m-1}$ leads instead to the constrained spline of Theorem~\ref{Theo:MMSE}, which annihilates these polynomials.

Our second variant of Theorem \ref{Theo:GenGaussMarkov} finalizes this functional picture by providing the  corresponding MMSE (or Bayesian) estimator in full generality. The latter also fits the description of a generalized spline, but of a constrained type that fulfills prescribed boundary conditions. For $\V p=\{0\}$ (trivial null space), we recover Theorem \ref{Theo:GenGaussMarkov}, with the MMSE and BLU estimators being equivalent.

\begin{theorem}[MMSE estimation of a GGP] 
\label{Theo:MMSE}
We consider the following setting.
\begin{itemize}
\item The signal $g=\{g(\varphi)\}_{\varphi \in \Spc S} \in \Spc S'$ to be recovered is a
realization of a generalized Gaussian process $\gRand g\sim \Spc N(0,\Op C_{\Lop,\V \phi})$
in $\Spc S'$ whose covariance operator $\Op C_{\Lop,\V \phi}=\Op A_{\V \phi}$ is the Gram
operator \eqref{III-Eq:Aphi} of a $\V p$-admissible operator $\Lop$; by
Theorem \ref{Theo:GGPcondi}, the action of $\gRand g$ extends to test functions
$\varphi \in \Spc H'_\Lop$.

\item The observed data are $\M y=\V \nu(g) + \V \epsilon \in \R^M$,
with $\V \nu=(\nu_1, \dots,\nu_M)\in (\Spc H'_\Lop)^M$ and $\V \epsilon$ a realization of a white Gaussian noise  with variance $\sigma^2$.
\end{itemize}
Then, for any $\varphi\in \Spc H'_\Lop\supset \Spc S$, the MMSE estimator of
$g(\varphi)$ given $\M y=(y_m)$
is 
\begin{align}
\label{Eq:MMSEsol00}
g_{{\rm MMSE}}(\varphi|\M y)=\Exp\{\gRand g( \varphi)|\M y\}%
&=\tilde g_{\lambda}( \varphi)=\langle \tilde g_{\lambda},\varphi\rangle\quad \mbox{ with }\quad  \lambda=\sigma^2,
\end{align}
where
$\tilde g_{\lambda}\in \Spc H_\Lop$ is the constrained generalized spline specified by
\begin{align}
\tilde g_\lambda&=\arg \min_{g \in \Spc H_\Lop} \left(\sum_{m=1}^M \big|y_m - \langle \nu_m,g\rangle\big|^2 + \lambda \|\Lop g\|^2_{L_2}\right) \mbox { s.t. } \V \phi(g)=0\label{Eq:GSmoothingSpline22}
\\
&=\sum_{m=1}^M a_m \Op A_\V \phi\{\nu_m\}\quad \mbox{ with } \quad (a_m)=(\widetilde{\M G} + \lambda \M I)^{-1} \M y, 
\label{Eq:GSmoothingSpline3}
\end{align}
where  $\widetilde{\M G} \in \R^{M \times M}$ with $[\widetilde{\M G}]_{m,n}=\langle \Op A_\V \phi \nu_m, \nu_n\rangle$.
\end{theorem}

\begin{proof}
Let $\widetilde{\varphi}=\Op P_{\V \phi}\{\varphi\}$ and $\widetilde{\nu}_m=\Op P_{\V \phi}\{\nu_m\}$
denote the projections onto $\Spc H_{\Op A}$. By Theorem \ref{Theo:GGPcondi},
$\gRand g(\varphi)=\gRand g(\widetilde{\varphi})$ and $\gRand g(\nu_m)=\gRand g(\widetilde{\nu}_m)$
almost surely, so that the measurement model is almost surely
$\M y=\widetilde{\V \nu}(g)+\V \epsilon$ and the quantity to estimate is
$\gRand g(\widetilde{\varphi})$. We now read this model on the dual pair
$(\Spc H_{\Lop,\V \phi},\Spc H_{\Op A})$, whose Riesz map is $\Op A_{\V \phi}$
\big(see \eqref{Eq:RieszPairingIII}\big). On this pair, $\Lop: \Spc H_{\Lop,\V \phi} \to L_2$
is invertible with inverse $\Op L^{-1}_{\V \phi}$, the norm is
$\|\tilde g\|_{\Spc H_{\Lop,\V \phi}}=\|\Lop \tilde g\|_{L_2}$ by \eqref{III-Eq:InnerProd},
and the restriction of $\gRand g$ to $\Spc H_{\Op A}$ is the non-degenerate GGP with
covariance $\Op A_{\V \phi}=(\Lop^\ast\Lop)^{-1}$, by \eqref{III-Eq:Factorization} and
\eqref{Eq:Cov2}. The derivation of Theorem \ref{Theo:GenGaussMarkov} uses nothing but these
three ingredients, and it therefore applies verbatim: with $\lambda=\sigma^2$,
\begin{align*}
\Exp\{\gRand g(\varphi)|\M y\}&=\Exp\{\gRand g(\widetilde{\varphi})|\M y\}
=\tilde g_\lambda(\widetilde{\varphi}),\\
\tilde g_\lambda&=\sum_{m=1}^M a_m \Op A_{\V \phi}\{\widetilde{\nu}_m\}
\quad \mbox{ with } \quad \M a=(\widetilde{\M G}+\lambda \M I)^{-1}\M y,
\end{align*}
where $\tilde g_\lambda$ is the minimizer of
$\sum_{m=1}^M |y_m-\langle \widetilde{\nu}_m,\tilde g\rangle|^2+\lambda\|\Lop \tilde g\|^2_{L_2}$
over $\Spc H_{\Lop,\V \phi}=\{g \in \Spc H_\Lop: \V \phi(g)=\V 0\}$. It remains to remove the
tildes. Since $\Spc N_{\V \phi}$ annihilates $\Spc H_{\Lop,\V \phi}$, we have that
$\langle \widetilde{\nu}_m,\tilde g\rangle=\langle \nu_m,\tilde g\rangle$ and
$\tilde g_\lambda(\widetilde{\varphi})=\tilde g_\lambda(\varphi)$, which identifies the
minimization problem with \eqref{Eq:GSmoothingSpline22}. Likewise,
$\Op A_{\V \phi}\{\widetilde{\nu}_m\}=\Op A_{\V \phi}\{\nu_m\}$ and
$\langle \Op A_{\V \phi}\{\widetilde{\nu}_m\},\widetilde{\nu}_n\rangle
=\langle \Op A_{\V \phi}\{\nu_m\},\nu_n\rangle=[\widetilde{\M G}]_{m,n}$, because
$\Op A_{\V \phi}=\Op A_{\V \phi}\Op P_{\V \phi}$ by \eqref{III-Eq:Aphi}, which gives
\eqref{Eq:GSmoothingSpline3}. Let us note that, for $\lambda>0$, no linear independence of
the $\widetilde{\nu}_m$ is needed, since $\widetilde{\M G}+\lambda \M I$ is always invertible.
\end{proof}

To the best of our knowledge, Theorem \ref{Theo:MMSE} is new at this level of generality,
precisely because its formulation requires the proper definition of the corresponding GGP.

\begin{remark}[Exact fit]
\label{Rem:ExactFit}
Theorem \ref{Theo:MMSE} remains valid in the noise-free scenario $\sigma^2=\lambda=0$,
provided that the Gram matrix $\widetilde{\M G}$ is invertible. Indeed, the conditioning
argument of the proof only requires the invertibility of
$\M C_{\widetilde{\VRand Y}}=\widetilde{\M G}+\sigma^2\M I$, while the coefficients become
$\M a=\widetilde{\M G}^{-1}\M y$, so that the estimator fits the data exactly, with
$\langle \nu_m,\tilde g_0\rangle=[\widetilde{\M G}\M a]_m=y_m$. The corresponding variational
characterization is the constrained interpolation problem
\begin{align*}
\tilde g_0=\arg\min_{g\in\Spc H_\Lop}\|\Lop g\|_{L_2}\quad\mbox{s.t.}\quad
\V\nu(g)=\M y \ \mbox{ and } \ \V\phi(g)=\V 0,
\end{align*}
which follows from \cite[Theorem~\ref{II-Theo:MinNorm}]{Unser2026b} applied to
$\Spc H_{\Lop,\V \phi}$, whose Riesz map is $\Op A_{\V \phi}$.
\end{remark}

For the limit case $\lambda\to0$, we can formulate a simpler version of
Theorem \ref{Theo:MMSE} that fills the gap that was stated in the introduction by showing that {\em every generalized spline interpolant is the MMSE estimator of a corresponding GGP}. 
To that end, we impose the condition $\Spc N_{\V \phi} \subset \Spc N_{\V \nu}={\rm span}\{\nu_m\}_{m=1}^M$, which means that the observations of the process must be rich enough to span the boundary conditions $\V \phi(\gRand G)=\V 0$.
\begin{corollary}[MMSE interpolation of a GGP]
\label{Theo:Unser}
Let $\Rand Y_m=\gRand G(\nu_m), m=1,\dots,M$ be a finite collection of observations of a generalized Gaussian process $\gRand G\sim \Spc N(0,\Op A_{\V \phi})$ in $\Spc S'$ that is whitened by a $\V p$-admissible operator $\Lop$, where the functionals $\nu_1,\dots,\nu_M \in \Spc H'_\Lop$ are linearly independent. Then, under the condition that $\Spc N_{\V \phi} \subset {\rm span}\{ \nu_m\}_{m=1}^{M}$, the {\em MMSE interpolator} of $\gRand G$ given the observed data
$\VRand Y=\M y$ is the generalized spline interpolant 
\begin{align}
g_0=\arg \min_{g \in \Spc H_\Lop} \|\Lop g\|_{L_2}
\quad \mbox{s.t.}\quad \V \nu(g_0)=\M y.
\end{align}
It is such that $\Exp\{\gRand G(\varphi)|\M y\}=g_0(\varphi)$ for any $\varphi \in \Spc H'_\Lop\supset \Spc S$ and almost every realization $\M y$ of the data.
\end{corollary}
\begin{proof}
The hypothesis provides a matrix $\widetilde{\M P} \in \R^{M \times N_0}$ such that
$\widetilde{\M P}^\Tr \V \nu=\V \phi$; since
$\widetilde{\M P}^\Tr \M P=\V \phi(\V p)=\M I_{N_0}$, the matrix $\M P=\V \nu(\V p)$ has
rank $N_0$ and \cite[Proposition~\ref{II-Prop:DirectSum}]{Unser2026b} applies with the choice
$\widetilde{\V \phi}=\V \phi$. It provides the change of coordinates
$\V \nu=\M Q\widetilde{\V \nu}+\M P \V \phi$, in which the data are equivalent to the pair
$\big(\widetilde{\M y},\widetilde{\M P}^\Tr \M y\big)
=\big(\widetilde{\M Q}^\Tr\M y,\V \phi(\gRand G)\big)$, with
$\widetilde{\V \nu}=\widetilde{\M Q}^\Tr\V \nu \in \Spc H_{\Op A}^{M-N_0}$. The second
block vanishes almost surely (Theorem \ref{Theo:GGPcondi}), so that it carries no
information---the $\sigma$-algebras generated by $\VRand Y$ and by
$\widetilde{\VRand Y}=\widetilde{\V \nu}(\gRand G)$ coincide---while, on the deterministic
side, it turns the interpolation constraints $\V \nu(g)=\M y$ into
$\widetilde{\V \nu}(g)=\widetilde{\M y}$ together with the boundary conditions
$\V \phi(g)=\V 0$. In particular, the feasible set of the interpolation problem is almost
surely contained in $\Spc H_{\Lop,\V \phi}$, so that $g_0$ coincides with the {\em
constrained} interpolant of the reduced data $\widetilde{\M y}$.

It then suffices to apply Theorem \ref{Theo:MMSE} with an exact fit
(Remark \ref{Rem:ExactFit}) to the reduced, noise-free measurement model
$\widetilde{\VRand Y}=\widetilde{\V \nu}(\gRand G)$. The corresponding Gram matrix
$\widetilde{\M G}=\big[\langle \Op A_{\V \phi}\{\widetilde{\nu}_m\},\widetilde{\nu}_n\rangle\big]
\in \R^{(M-N_0)\times(M-N_0)}$ is invertible, because the $\widetilde{\nu}_m$ are linearly
independent \big($\widetilde{\M Q}$ has full column rank and the $\nu_m$ are linearly
independent\big) and $\Op A_{\V \phi}$ is the Riesz map of $\Spc H_{\Op A}$. Hence
$\Exp\{\gRand G(\varphi)\,|\,\widetilde{\VRand Y}=\widetilde{\M y}\}=\tilde g_0(\varphi)$,
where $\tilde g_0$ is the constrained interpolant of $\widetilde{\M y}$; that is,
$\tilde g_0=g_0$.
\end{proof}

\section{Conclusion}

We conclude the paper by emphasizing the crucial role of the regularization/whitening operator $\Lop$ and of its native space $\Spc H_\Lop$. As for the deterministic part of the story, the generalized forms of $\Lop$-splines are the solutions of regularized linear inverse problems in $\Spc H_\Lop$. Their stochastic counterparts are precisely the GGPs whitened by $\Lop$. (The delicate part in this statement is the necessity to invert $\Lop$, which requires the specification of boundary conditions.) The pleasing outcome is the mathematical congruence between the deterministic and stochastic worlds. The present framework provides a rigorous support to statements such as ``fractional splines are optimal estimators for fractional Brownian motion.'' Likewise, we can now identify the polynomial splines of degree $2n-1$ as the optimal estimators of Gaussian processes, which, up to some boundary conditions (the technical part of the story), are the $n$-fold integration of a white noise, or, for more classic folks,  the $(n-1)$-fold integration of a Brownian motion.
With the help of abstraction and nuclear spaces, we ultimately end up with an equivalence in the continuum that is as strong as in finite dimensions.

\appendix

\section{Higher-Order fBm and Fractional Splines}
\label{App:fBm}
Section \ref{Sec:GaussianSDE} presented fractional Brownian motion (fBm) as the
prototypical generalized Gaussian process with a one-dimensional null space. We now carry
out the same construction for an arbitrary order, which yields the family of
{\em higher-order fractional Brownian motions} \cite{Perrin2001}, and we make the three
objects around which the theory revolves fully explicit: the conditioned inverse
$\Op L^{-1}_{\V \phi}$, the covariance operator $\Op A_{\V \phi}$, and the covariance
function of the process. The biorthogonal system is the one that
\cite[Table~\ref{II-Tab:Derivatives}]{Unser2026b} attaches to fractional splines, so that the boundary
conditions are the vanishing of the function and of its first $(N_0-1)$ derivatives at the
origin.

\subsection{Operator, Boundary Conditions, and Admissible Range}
\label{Sec:fBmSetting}
The relevant whitening operator is the fractional derivative $\Lop=\Dop^{\gamma}: \Spc S(\R) \to \Spc S'(\R)$
of order $\gamma>\tfrac{1}{2}$, whose Fourier symbol is
$\widehat{L}(\omega)=(\jj \omega)^{\gamma}$. We also set
\begin{align}
\label{Eq:HurstGamma}
H=\gamma-\tfrac{1}{2},
\end{align}
which is the Hurst-type exponent of the resulting process. The biorthogonal system of
\cite[Table~\ref{II-Tab:Derivatives}]{Unser2026b} is
\begin{align}
\label{Eq:fBmBiortho}
p_n(x)=\frac{x^{n-1}}{(n-1)!}, \qquad \phi_n=\delta^{(n-1)}, \qquad n=1,\dots,N_0,
\end{align}
where we adopt the convention $\langle \delta^{(n-1)},f\rangle=f^{(n-1)}(0)$, which is the
one that makes $(\V p,\V \phi)$ biorthogonal with
$\langle p_m,\phi_n\rangle=\delta_{m,n}$. The null space
$\Spc N_{\V p}={\rm span}\{1,x,\dots,x^{N_0-1}\}$ is the space of polynomials of degree
less than $N_0$, while the boundary conditions that specify the conditioned inverse read
\begin{align}
\label{Eq:fBmBC}
\V \phi(g)=\V 0 \quad \Longleftrightarrow \quad g(0)=g'(0)=\cdots=g^{(N_0-1)}(0)=0.
\end{align}
These are the standard homogeneous initial conditions for the function and its
derivatives, and they reduce to $g(0)=0$ for $N_0=1$, which is Mandelbrot's normalization
of fBm.

\begin{remark}[Admissible range]
\label{Rem:fBmRange}
The pair $(\gamma,N_0)$ cannot be chosen independently. The quadratic form
\eqref{Eq:fBmQuadratic} below is finite for every $\varphi \in \Spc S(\R)$ if and only if
\begin{align}
\label{Eq:fBmRange}
N_0-\tfrac{1}{2}<\gamma<N_0+\tfrac{1}{2},
\end{align}
which identifies $N_0=\lfloor \gamma+\tfrac{1}{2}\rfloor$ as the integer that is nearest to
$\gamma$, in agreement with the general rule $\Spc N_{\Op B}=\Pi_{k_\gamma}$ with
$k_\gamma=\lfloor \gamma-\tfrac{d}{2}\rfloor$ of \cite{Unser2026a} and with the entry for $\Dop^\alpha$ of
\cite[Table~\ref{II-Tab:Derivatives}]{Unser2026b}. The lower bound is what neutralizes the
singularity of $|\omega|^{-2\gamma}$ at the origin, while the upper bound is what keeps the
correction term integrable at infinity; equivalently, the upper bound is the condition under
which the monomial $x^{N_0-1}$ belongs to the attainable null space of $\Op D^{\gamma\ast}\Op D^{\gamma}$. The two boundary cases
$\gamma \in \N+\tfrac{1}{2}$ are excluded, and they are precisely the ones for which the
constant $c_\gamma$ of \eqref{Eq:rhogamma} below becomes singular and the Green's function
acquires a logarithmic factor. For $\gamma=m \in \N$, \eqref{Eq:fBmRange} yields $N_0=m$,
which is the classical polynomial-spline setting; for the classical range
$H \in (0,1)$ of fBm, it yields $N_0=1$ for every $H$.
\end{remark}

\subsection{The Conditioned Inverse and its Fourier Characterization}
\label{Sec:fBmInverse}
The Green's operator of $\Dop^\gamma$ is the Riemann--Liouville fractional integral
$\Op D^{-\gamma}: w \mapsto \rho_\Lop \ast w$ with $\rho_\Lop(x)=\frac{x_+^{\gamma-1}}{\Gamma(\gamma)}$
\cite[Table~\ref{II-Tab:Derivatives}]{Unser2026b}, while the conditioning is carried by the projector
$\Op P^\ast_{\V \phi}=\Identity-\Op R_{\V p}\Op R_{\V \phi}$ of Section \ref{Sec:NativeRecap}, whose action on a
sufficiently smooth function is the subtraction of its Taylor polynomial at the origin:
\begin{align}
\label{Eq:TaylorProj}
\Op P^\ast_{\V \phi}\{f\}(x)=f(x)-\sum_{n=1}^{N_0} \frac{x^{n-1}}{(n-1)!}\, f^{(n-1)}(0).
\end{align}
The conditioned inverse $\Op L^{-1}_{\V \phi}=\Op P^\ast_{\V \phi}\Op D^{-\gamma}$ is
therefore an integral operator with the (non-convolutional) kernel
\begin{align}
\label{Eq:fBmKernel}
\rho_{\gamma,\V \phi}(t,\tau)=\frac{(t-\tau)_+^{\gamma-1}}{\Gamma(\gamma)}
-\sum_{n=1}^{N_0} \frac{t^{n-1}}{(n-1)!}\,\frac{(-\tau)_+^{\gamma-n}}{\Gamma(\gamma-n+1)},
\end{align}
which is obtained by applying \eqref{Eq:TaylorProj} to $t \mapsto \rho_\Lop(t-\tau)$, and
which is the instance of the general formula \eqref{Eq:kernelinverseg} of
Theorem \ref{Theo:ReproKernel} that corresponds to $q_n(\tau)=\frac{(-\tau)_+^{\gamma-n}}{\Gamma(\gamma-n+1)}$. For
$N_0=1$, \eqref{Eq:fBmKernel} reduces to
$\frac{1}{\Gamma(\gamma)}\big((t-\tau)_+^{\gamma-1}-(-\tau)_+^{\gamma-1}\big)$, which is
the integral kernel of Mandelbrot's representation \eqref{Eq:fBMH}. The higher-order
kernels are its natural extension, with one correction term per boundary condition.

The conditioning breaks the shift-invariance, so that $\Op L^{-1}_{\V \phi}$ has no
frequency response of its own. Its adjoint, however, has a transparent Fourier
description. Let $\hat \varphi=\Fourier\{\varphi\}$ and let
\begin{align}
\label{Eq:TaylorFourier}
\hat \varphi_{N_0}(\omega)=\sum_{n=1}^{N_0} \frac{\hat \varphi^{(n-1)}(0)}{(n-1)!}\,\omega^{n-1}
\end{align}
be the Taylor polynomial of $\hat \varphi$ of degree $(N_0-1)$ at the origin. Since
$\Fourier\{\phi_n\}(\omega)=(-\jj \omega)^{n-1}$ and
$\langle p_n,\varphi\rangle=\frac{\jj^{\,n-1}}{(n-1)!}\hat \varphi^{(n-1)}(0)$, the
projector $\Op P_{\V \phi}=\Identity-\Op R_{\V \phi}\Op R_{\V p}$ acts in the Fourier
domain as
\begin{align}
\label{Eq:ProjFourier}
\Fourier\{\Op P_{\V \phi}\varphi\}(\omega)=\hat \varphi(\omega)-\hat \varphi_{N_0}(\omega),
\end{align}
which is the exact Fourier-domain counterpart of \eqref{Eq:TaylorProj}. Together with
$\Op L^{-1\ast}_{\V \phi}=\Op D^{-\gamma\ast}\Op P_{\V \phi}$, this yields the announced
characterization
\begin{align}
\label{Eq:fBmInverseFourier}
\Fourier\big\{\Op L^{-1\ast}_{\V \phi}\varphi\big\}(\omega)
=\frac{\hat \varphi(\omega)-\hat \varphi_{N_0}(\omega)}{(-\jj \omega)^{\gamma}},
\qquad (-\jj \omega)^{\gamma}=\overline{(\jj \omega)^{\gamma}}.
\end{align}
The interpretation is that the conditioning removes the first $N_0$ Taylor coefficients of
$\hat \varphi$ at the origin, which is exactly what is needed for the division by
$(\jj \omega)^{\gamma}$ to remain well-defined. In signal-processing terms, the operator
$\Op P_{\V \phi}$ endows the test function with $N_0$ vanishing moments.

\subsection{The Covariance Operator}
\label{Sec:fBmCovOp}
By combining \eqref{Eq:fBmInverseFourier} with the factorization
$\Op A_{\V \phi}=\Op L^{-1}_{\V \phi}\Op L^{-1\ast}_{\V \phi}$ \big(see \eqref{III-Eq:Factorization}\big),
we obtain the Fourier-domain description of the covariance operator as the sesquilinear
form
\begin{align}
\label{Eq:fBmQuadratic}
\langle \Op A_{\V \phi} \varphi,\psi\rangle
=\frac{1}{2\pi}\int_\R \frac{\big(\hat \varphi(\omega)-\hat \varphi_{N_0}(\omega)\big)
\overline{\big(\hat \psi(\omega)-\hat \psi_{N_0}(\omega)\big)}}{|\omega|^{2\gamma}}\,\dint \omega,
\end{align}
which is finite for all $\varphi,\psi \in \Spc S(\R)$ under the condition
\eqref{Eq:fBmRange}. The integrand is non-negative for $\psi=\varphi$, and it vanishes
identically exactly on $\Spc N_{\V \phi}={\rm span}\{\delta^{(n-1)}\}_{n=1}^{N_0}$, whose
members have the polynomial Fourier transforms $\Fourier\{\phi_n\}(\omega)=(-\jj \omega)^{n-1}$:
this is the degeneracy of Theorem \ref{Theo:GGPcondi} in concrete form. For
$\varphi,\psi \in \Spc H_{\Op A}$, the corrections disappear, because
$\langle p_n,\varphi\rangle=0$ is equivalent to $\hat \varphi^{(n-1)}(0)=0$, and
\eqref{Eq:fBmQuadratic} reduces to
$\langle \Op A \varphi,\psi\rangle=\frac{1}{2\pi}\int_\R \frac{\hat \varphi \overline{\hat \psi}}{|\omega|^{2\gamma}}\dint \omega$,
which is the concrete form of \eqref{Eq:EqinnerProdHA}.

\subsection{The Covariance Function}
\label{Sec:fBmCovFunction}
The process $\gRand B_\gamma=\Op L^{-1}_{\V \phi}\{\gRand W\} \sim \Spc N(0,\Op C_{\Lop,\V \phi})$
is the solution of the fractional stochastic differential equation
\begin{align}
\label{Eq:fBmSDE}
\Dop^{\gamma}\{\gRand B_\gamma\}=\gRand W \qquad \mbox{ subject to } \qquad
\gRand B_\gamma^{(n-1)}(0)=0,\quad n=1,\dots,N_0,
\end{align}
which is the concrete instance of \eqref{Eq:SDE}. For $\gamma>\tfrac{1}{2}$, the native
space $\Spc H_{\Dop^\gamma}$ is a reproducing-kernel Hilbert space of continuous functions,
with $\delta_t \in \Spc H'_{\Dop^\gamma}$ for every $t \in \R$; by
Item \ref{Item:RKHS} of Theorem \ref{Theo:ReproKernel}, this is equivalent to the
finiteness of the diagonal value $a_{\V \phi}(t,t)=\sigma^2_\gamma\,|t|^{2H}$ computed in
\eqref{Eq:fBmVariance} below. The covariance function
\begin{align}
\label{Eq:fBmCovDef}
r_\gamma(t,\tau)=\Exp\big\{\gRand B_\gamma(t)\gRand B_\gamma(\tau)\big\}
=\langle \Op A_{\V \phi}\{\delta_\tau\},\delta_t\rangle
\end{align}
is well-defined. Let
\begin{align}
\label{Eq:rhogamma}
\rho_\gamma(x)=c_\gamma\,|x|^{2\gamma-1}, \qquad
c_\gamma=\frac{-1}{2\sin\big(\tfrac{\pi}{2}(2\gamma-1)\big)\Gamma(2\gamma)}
=\frac{1}{2\cos(\pi \gamma)\Gamma(2\gamma)},
\end{align}
be the Green's function of $\Op D^{\gamma\ast}\Op D^{\gamma}$ of
\cite[Table~\ref{II-Tab:Derivatives}]{Unser2026b}, whose derivatives are
$\rho^{(k)}_\gamma(x)=c_\gamma\frac{\Gamma(2\gamma)}{\Gamma(2\gamma-k)}|x|^{2\gamma-1-k}({\rm sgn}\,x)^k$.
\begin{proposition}[Covariance of the higher-order fBm]
\label{Theo:fBmCov}
Under the condition \eqref{Eq:fBmRange}, the covariance function of the process
\eqref{Eq:fBmSDE} is
\begin{align}
\label{Eq:fBmCovGeneral}
r_\gamma(t,\tau)=\rho_\gamma(t-\tau)
-\sum_{n=1}^{N_0}\frac{(-t)^{n-1}}{(n-1)!}\rho^{(n-1)}_\gamma(\tau)
-\sum_{n=1}^{N_0}\frac{(-\tau)^{n-1}}{(n-1)!}\rho^{(n-1)}_\gamma(t).
\end{align}
It is the double Taylor remainder of $\rho_\gamma(t-\tau)$, in the sense that it is
obtained by subtracting from $\rho_\gamma(t-\tau)$ its Taylor polynomials of degree
$(N_0-1)$ at $t=0$ and at $\tau=0$.
\end{proposition}
\begin{proof}
This is the general formula \eqref{Eq:kernelrepro} of Theorem \ref{Theo:ReproKernel}
specialized to $\Lop=\Dop^\gamma$, whose double sum vanishes here. We give the direct
computation. The identity $\Op A_{\V \phi}=\Op P^\ast_{\V \phi}\Op A\Op P_{\V \phi}$ of \eqref{III-Eq:Aphi} gives
$r_\gamma(t,\tau)=\langle \Op A\,\Op P_{\V \phi}\{\delta_\tau\},\Op P_{\V \phi}\{\delta_t\}\rangle$
with $\Op P_{\V \phi}\{\delta_t\}=\delta_t-\sum_{n=1}^{N_0}p_n(t)\,\phi_n$. The expansion
of the bilinear form produces four groups of terms, which we evaluate by means of
$\langle \Op A \delta_\tau,\delta_t\rangle=\rho_\gamma(t-\tau)$,
$\langle \Op A \delta_\tau,\phi_n\rangle=\rho^{(n-1)}_\gamma(-\tau)$,
$\langle \Op A \phi_m,\delta_t\rangle=(-1)^{m-1}\rho^{(m-1)}_\gamma(t)$, and
$\langle \Op A \phi_m,\phi_n\rangle=(-1)^{m-1}\rho^{(m+n-2)}_\gamma(0)$. The last group
vanishes, because $\rho_\gamma$ is $(m+n-2)$ times differentiable at the origin with a
vanishing derivative there whenever $m+n<2\gamma+1$, a condition that
\eqref{Eq:fBmRange} guarantees for all $m,n \le N_0$. In the notation of
Theorem \ref{Theo:ReproKernel}, this is the statement that the Gram coefficients
\eqref{Eq:rmn} all vanish, so that $r_{m,n}=0$. The evenness of $\rho_\gamma$ gives
$\rho^{(n-1)}_\gamma(-\tau)=(-1)^{n-1}\rho^{(n-1)}_\gamma(\tau)$, which turns the two
middle groups into the two sums of \eqref{Eq:fBmCovGeneral}.
\end{proof}
One verifies directly on \eqref{Eq:fBmCovGeneral} that
$\partial^{n-1}_t r_\gamma(t,\tau)|_{t=0}=0$ for $n=1,\dots,N_0$, which is the covariance
transcription of the boundary conditions \eqref{Eq:fBmBC}.

For $N_0=1$, that is, for $\gamma \in (\tfrac{1}{2},\tfrac{3}{2})$ and $H \in (0,1)$,
formula \eqref{Eq:fBmCovGeneral} reduces to the classical covariance of fBm,
\begin{align}
\label{Eq:fBmCovN1}
r_\gamma(t,\tau)=\frac{\sigma^2_H}{2}\big(|t|^{2H}+|\tau|^{2H}-|t-\tau|^{2H}\big),
\qquad
\sigma^2_H=-2c_\gamma=\frac{1}{\sin(\pi H)\Gamma(2H+1)},
\end{align}
which also fixes the normalization constant that is implicit in \eqref{Eq:CovfBm}: it is
the one that corresponds to the choice
$\rho_\Lop(x)=\frac{x_+^{\gamma-1}}{\Gamma(\gamma)}$ of \cite[Table~\ref{II-Tab:Derivatives}]{Unser2026b}. For
$H=\tfrac{1}{2}$, one gets $\sigma^2_H=1$ and $r_1(t,\tau)=\min(t,\tau)$ for
$t,\tau \ge 0$, which is Brownian motion. For $N_0=2$, that is, for
$\gamma \in (\tfrac{3}{2},\tfrac{5}{2})$ and $H \in (1,2)$, we get
\begin{align}
\label{Eq:fBmCovN2}
r_\gamma(t,\tau)=c_\gamma\Big(|t-\tau|^{2H}-|t|^{2H}-|\tau|^{2H}
+2H\big(t\,|\tau|^{2H-1}{\rm sgn}(\tau)+\tau\,|t|^{2H-1}{\rm sgn}(t)\big)\Big),
\end{align}
where the two additional terms are the ones that enforce
$\gRand B_\gamma'(0)=0$.

Whatever the order, the process remains self-similar of exponent $H$: the substitution
$\tau=t$ in \eqref{Eq:fBmCovGeneral} gives
\begin{align}
\label{Eq:fBmVariance}
\Exp\big\{|\gRand B_\gamma(t)|^2\big\}=\sigma^2_\gamma\,|t|^{2H},
\qquad
\sigma^2_\gamma=-2c_\gamma\,\Gamma(2\gamma)\sum_{n=1}^{N_0}\frac{(-1)^{n-1}}{(n-1)!\,\Gamma(2\gamma-n+1)},
\end{align}
which is consistent with $\sigma^2_\gamma=\sigma^2_H$ for $N_0=1$.

Formula \eqref{Eq:fBmCovGeneral} coincides with the covariance that Perrin et al.\
obtained for their $n$th-order extension of fractional Brownian motion \cite{Perrin2001},
with $n=N_0$ and with the same normalization: the vanishing of the process and of its
first $(N_0-1)$ derivatives at the origin. The transfer principle of \cite{Sottinen2019}
relates these processes to the classical fBm.

\subsection{The Associated Splines}
\label{Sec:fBmSplines}
Theorem \ref{Theo:MMSE} now applies verbatim. Given the noisy samples
$y_m=\gRand B_\gamma(t_m)+\epsilon_m$ at $M$ distinct locations $t_1,\dots,t_M$, the MMSE
estimator of $\gRand B_\gamma$ is the constrained fractional spline
\begin{align}
\label{Eq:fBmSpline}
\tilde g_\lambda=\sum_{m=1}^M a_m\, r_\gamma(\cdot,t_m),
\qquad \M a=(\widetilde{\M G}+\lambda \M I)^{-1}\M y, \qquad \lambda=\sigma^2,
\end{align}
with $[\widetilde{\M G}]_{m,n}=r_\gamma(t_m,t_n)$. The basis functions are the covariance
function itself, in accordance with the congruence of
Proposition \ref{Theo:Congruence}. For $\gamma=m \in \N$, one has $N_0=m$,
$c_\gamma=\frac{(-1)^m}{2(2m-1)!}$, and \eqref{Eq:fBmSpline} is the polynomial smoothing
spline of degree $(2m-1)$ that satisfies the boundary conditions
\eqref{Eq:fBmBC}; for non-integer $\gamma$, it is the fractional spline of
\cite{Unser2000,Unser2000a}. This gives a precise meaning, and a proof, to the statement
that fractional splines are the optimal estimators of fractional Brownian motion
\cite{Blu2007a,Tirosh2006}.

The interpolation counterpart (Corollary \ref{Theo:Unser}) is more demanding, since its
hypothesis $\Spc N_{\V \phi} \subset \Spc N_{\V \nu}$ requires the boundary functionals to
be reproducible from the observations. For $N_0=1$, it suffices to include the sample at
$t=0$, whose value $\gRand B_\gamma(0)=0$ is deterministic. For $N_0 \ge 2$, the
derivatives $\delta^{(n-1)}$ cannot be reproduced by finitely many point evaluations, so
that the samples must be complemented by the derivative measurements
$\gRand B^{(n-1)}_\gamma(0)=0$ for $n=2,\dots,N_0$. Once this is done, every fractional
spline interpolant of the augmented data is an MMSE interpolator.

\subsection*{Acknowledgment}
The author used Claude (Anthropic) as an interactive assistant during the preparation of
this paper, for proofreading, for checking and simplifying proofs, and for the numerical
verification of the formulas of the appendix and Section \ref{Sec:fBf}. All results were verified by the author,
who takes full responsibility for the content.
This research was supported by the Swiss National Science Foundation (SNSF) under Grant 200020\_219356.

\bibliographystyle{plaininitials}
\bibliography{SplineOps}

\end{document}